\documentclass[review]{elsarticle}

\usepackage{lineno,hyperref}
\modulolinenumbers[5]

\usepackage{graphicx}
\usepackage{subfig}
\usepackage{amsmath}
\usepackage{lscape} 

\usepackage{nomencl} 
\usepackage{multicol}
\makenomenclature
\usepackage{etoolbox}
\renewcommand\nomgroup[1]{%
  \item[\bfseries
  \ifstrequal{#1}{A}{List of Symbols}{%
  \ifstrequal{#1}{B}{List of Subscripts}{}}%
]}

\journal{Journal of Computational and Applied Mathematics, }

\begin{document}

\begin{frontmatter}

\title{Efficient Radial Basis Function Mesh Deformation Methods for Aircraft Icing}

\author[mymainaddress,mysecondaryaddress]{Myles Morelli\corref{mycorrespondingauthor}}
\cortext[mycorrespondingauthor]{Corresponding author, email address: mylescarlo.morelli@polimi.it}

\author[mymainaddress]{Tommaso Bellosta}

\author[mymainaddress]{Alberto Guardone}

\address[mymainaddress]{Department of Aerospace Science and Technology, Politecnico di Milano, Italy}
\address[mysecondaryaddress]{CFD Laboratory, School of Engineering, University of Glasgow, United Kingdom}

\begin{abstract}
This paper presents an evaluation of efficient radial basis function mesh deformation for complex iced geometries. Given the high computational cost of mesh deformation, state-of-the-art radial basis function techniques are used for data reduction. The principle procedures adopted are multi-level greedy surface point selection and volume point reduction. The multi-level greedy surface point selection reduces the control point list to increase the efficiency of the interpolation operation and the volume point reduction improves the computational cost of the volume mesh update operation which is important for large data sets. The study demonstrates the capabilities of radial basis function mesh deformation in both two and three-dimensions. Furthermore, it compares localised ice deformation to more standardized test cases with global deformation. The convergence history of the multi-level greedy point selection is assessed in terms of number of control points and computational cost for all the test cases. The location of the selected control points near the ice accretion illustrates the effectiveness of the method for localised deformation. The results show that the radial basis function mesh deformation performs well for both the two and three-dimensional test cases. The reduction of the relative surface error for the three-dimensional test cases understandably requires a larger number of control points and thus results in a higher computational cost. Nevertheless, the data-reduction schemes presented in this work represent a significant improvement to standard radial basis function mesh deformation for three-dimensional aircraft icing tests with large data-sets.
\end{abstract}


\end{frontmatter}

\nomenclature[A, 01]{$b$}{semi-span} 
\nomenclature[A, 02]{$c$}{chord} 
\nomenclature[A, 03]{$D$}{support distance of the wall distance function}                
\nomenclature[A, 04]{$\textbf{E}$}{error vector of surface displacements}
\nomenclature[A, 05]{$f(\textbf{r})$}{function to be evaluated at position \textbf{r}}
\nomenclature[A, 06]{$k$}{volume reduction factor}
\nomenclature[A, 07]{$N$}{number of points}
\nomenclature[A, 08]{$\textbf{r}$}{radial distance}
\nomenclature[A, 09]{$R$}{support radius of the radial basis function}
\nomenclature[A, 10]{$\Delta S$}{surface displacement}
\nomenclature[A, 11]{$\Delta V$}{volume displacement}
\nomenclature[A, 12]{$\Delta \textbf{X}$}{x-displacement vector of mesh coordinates}
\nomenclature[A, 13]{$\Delta \textbf{Y}$}{y-displacement vector of mesh coordinates}
\nomenclature[A, 14]{$\Delta \textbf{Z}$}{z-displacement vector of mesh coordinates}
\nomenclature[A, 15]{$\alpha$}{weight coefficients of the basis points}
\nomenclature[A, 16]{$\varepsilon$}{greedy tolerance}
\nomenclature[A, 17]{$\eta$}{ $ \left \| \textbf{r} - \textbf{r}_{i} \right \| / R $}
\nomenclature[A, 18]{$\xi$}{ $ d(\textbf{r})/D $}
\nomenclature[A, 19]{$\phi$}{radial basis function}
\nomenclature[A, 20]{$\Phi$}{universal basis matrix}
\nomenclature[A, 21]{$\psi$}{wall distance function}

\nomenclature[B]{$c$}{control points}
\nomenclature[B]{$l$}{multi-level}
\nomenclature[B]{$s$}{surface points}
\nomenclature[B]{$v$}{volume points}

\printnomenclature


\section{Introduction}
In-flight icing encounters can jeopardise the performance and handling qualities of aircraft and hence pose a serious threat to flight safety \cite{gent2000aircraft}. This threat to flight safety is brought forth as ice accretion can rapidly alter aerodynamic lifting surfaces such as wings and rotors during flight which are highly sensitive to geometric modifications. The ability to use computational techniques to simulate and model in-flight ice accretion has introduced alternative approaches to in-flight icing trails and experimental wind tunnel tests and is helping to further understand this safety-critical issue. However, the computational modelling of in-flight icing is not without its own set of challenges, and one challenge at the forefront is the issue of accounting for the moving ice boundary. Ice accretion can produce geometrically complex shapes, and with icing simulations being highly sensitive to the mesh quality, suitable mesh deformation schemes are required to maintain mesh orthogonality.

Well established techniques for mesh deformation such as the spring analogy \cite{batina1990unsteady, farhat1998torsional}, the linear elastic analogy \cite{baker1999dynamic, nielsen2002recent} and the interpolation method based on radial basis functions \cite{de2007mesh, rendall2009efficient, rendall2008unified, beckert2001multivariate} have been principally developed for moving boundary problem present in simulations such as fluid-structure interaction problems and aerodynamic shape optimization. The spring analogy, first developed by Batina \cite{batina1990unsteady} is one of the most widely used methods and models each edge of the mesh as a linear spring connected together at corresponding nodes. Significantly, Farhat et al. \cite{farhat1998torsional} further developed the spring analogy and introduced torsional stiffness to alleviate the mesh crossover problem. A concern with the spring analogy is that it is expensive to solve and cannot guarantee high mesh quality during large deformations. The linear elasticity analogy extends from the spring analogy and models each mesh element as an elastic solid. The linear elasticity analogy has very high mesh quality during large mesh deformations however this comes at the expense of an increase in the computational cost. Differently from large, smooth deformation present during shape optimisation or due to the effects of wing or airframe elasticity, geometry modifications due to ice accretion are usually local in nature and characterised by non-smooth shapes. Consequently, the accounting for the evolution of the iced surface boundary using conventional mesh deformation techniques is highlighted within the literature as being an issue \cite{tong2017three}.

Boer et al. first applied radial basis function interpolation to mesh deformation \cite{de2007mesh}. Radial basis functions mesh deformations techniques have the unique property that they do not require the grid connectivity meaning that even for three-dimensional problems they are relatively simple to implement. Their work observed the superior mesh quality of radial basis function interpolation when compared to the spring analogy \cite{de2007mesh}. Rendall and Allen \cite{rendall2008unified} then went on to show that radial basis function methods are robust and preserve high-quality mesh even during large deformations. One of the biggest deterrents to the radial basis function method, however, is that it is expensive for large scale problems when the number of surface points to be displaced becomes high. To address this issue and Rendall and Allen introduced data reduction schemes based upon greedy point selection algorithms to improve the efficiency of the interpolation operation \cite{rendall2009efficient, rendall2010reduced}. To further enhance this Wang et al. \cite{wang2015improved} introduced the concept of a multi-level subspace radial basis function interpolation method. To improve the volume update procedure Xie et al. \cite{xie2017efficient} developed a volume reduction scheme for large data sets. It is clear from the expedited literature that radial basis function mesh deformation is to play a prominent role in future mesh deformation advancements.  

This work seeks to explore how radial basis functions can be used with their robust mesh deformation properties for computational aircraft icing simulations during the mesh update procedure. Groth et al. \cite{groth2019rbf} first demonstrated the use of standard radial basis function mesh deformation techniques for icing simulations. Their work highlights the promise and potential of radial basis functions. However, significantly what it did not highlight was the critical issue of the efficiency of radial basis functions. Acknowledging the issue of the efficiency of the data-reduction schemes used in radial basis function mesh deformation is paramount and needs addressing. Within this work the most concurrent data reduction schemes are adopted and assessed for localised ice deformations. These methods are implemented within the open-source SU2 code which is used for solving the flow field for the icing simulations and the mesh deformation \cite{economon2016su2}. The ice accretion simulations are performed using the in-flight icing software suite developed by Politecnico di Milano called PoliMIce \cite{gori2015polimice}. 

The outline of the forthcoming paper is as follows: An introduction to radial basis functions and their use for mesh deformation is discussed in Section~\ref{sec:rbf}; The multi-level subspace data reduction schemes used in this work are described in Section~\ref{sec:data_reduction}; A volume reduction method for large data sets is discussed in Section~\ref{sec:volume}; The results of the methods applied to localised ice accretion specific test cases as well as more conventional global deformation benchmark test cases are presented in Section~\ref{sec:results}; The main talking points from this study are discussed in Section~\ref{sec:conclusion};

\section{Radial Basis Function Mesh Deformation}
\label{sec:rbf}
The term radial basis function refers to a series of functions whose values depends on their distance to a supporting position. In the most general of forms, radial basis functions can be written as, $ \phi (\textbf{r}, \textbf{r}_{i}) = \phi\left( \| \textbf{r} - \textbf{r}_{i} \|\right) $, where the distance corresponds to the radial basis centre, $ \textbf{r}_{i} $. This distance is frequently taken as the Euclidean distance, meaning it becomes the spatial distance between two nodes. 

An interpolation function, $f(\textbf{r})$ can be introduced as a method for describing the displacement of a set of nodes in space and can be approximated by a weighted sum of basis functions. However, the interpolation relies on the weight coefficients of the basis points, $\alpha$. The interpolation introduced by Ref.~\cite{rendall2008unified} takes the form

\begin{equation}
\label{eq:1}
    f(\textbf{r}) = \sum_{i=1}^{N} \alpha_{i} \phi \left ( \| \textbf{r} - \textbf{r}_{i}  \| \right ).
\end{equation}
The weight coefficients of the basis points described in equation~\ref{eq:1} can be obtained through the exact recovery of the known function values at the surface nodes. A prerequisite to this is knowing the surface node displacements a priori. The displacement of the surface points are contained within the vector, $ \Delta \textbf{X}$ as described by 

\begin{equation}
\label{eq:5}
    \Delta \textbf{X}_{s} = \left [ \Delta x_{s_1}, \, \Delta x_{s_2}, \, \dots , \, \Delta x_{N_s}  \right ]^T \, ,
\end{equation}
where the subscript ``$s$'' denotes the surface points. The $x$ components of the displacements are shown, however, the $y$ and $z$ components are analogous. This expression is reduced for the total number of surface points, $N_s$. The displacement in all $x$, $y$ and $z$ directions can be simplified to

\begin{equation}
    \Delta S = \Delta \textbf{X}_{s}\hat{\textbf{x}} + \Delta \textbf{Y}_{s}\hat{\textbf{y}} + \Delta \textbf{Z}_{s}\hat{\textbf{z}} \, .
\end{equation}
Similarly to the surface node displacements, the weight coefficients, $ \boldsymbol \alpha $, can be contained within a vector

\begin{equation}
    \boldsymbol \alpha_{x} = \left [\alpha_{x,s_1}, \, \alpha_{x,s_2}, \, \dots , \, \alpha_{x,N_s}  \right ]^T \, ,
\end{equation}
reiterating, only the $x$ components of the coefficients are shown, however, the $y$ and $z$ components are analogous. The universal basis matrix, $\boldsymbol \Phi$ can be constructed from the radial basis functions at each of the surface nodes and is consequently of the size of $N_s^2$. The universal basis matrix can then be shown in its compact form as

\begin{equation}
\label{eq:basisMatrix}
    \boldsymbol \Phi_{s_j, \, s_i} = \phi {\left \| r_{s_i} - r_{s_j} \right \|} \, .
\end{equation}
The coefficients can then be found by solving the following linear system:

\begin{equation}
\label{eq:9}
    \boldsymbol \Phi_{s, s} \boldsymbol \alpha = \Delta S \, ,
\end{equation}
the volume displacements, $\Delta V $, can finally be interpolated through the multiplication of the weight coefficients from equation~\ref{eq:9} and the newly formed volume node basis matrix, $\boldsymbol \Phi_{v, s}$, now of the size $N_v \times N_s$, as described by

\begin{equation}
\label{eq:10}
    \Delta V = \boldsymbol \Phi_{v, s} \boldsymbol \alpha \, ,
\end{equation}
where the subscript ``$v$'' represents the of volume points.

Multiple forms of radial basis functions exist within the literature which can be used for interpolating data sets and can be characterised into functions with global, local and, compact support. Functions with global support are always non-zero and grow with increasing distance from the radial basis function centre. Likewise, functions with local support are also always non-zero however decay with increasing distance from the radial basis function centre. Compact functions differ from global and local functions in that they decay to zero with increasing distance from the radial basis function centre. The choice of basis function is significant; global and local functions cover the entire interpolation space, leading to dense matrix systems which requires solving the linear system of a fully populated and ill-conditioned matrix. Compact functions are limited to a given support radius, $R$, resulting in sparse matrix systems which can be solved more easily however this sacrifices interpolation accuracy. 

With practical application in mind, functions with compact support were considered within this work such that, $\phi \left ( \left \| \textbf{r} - \textbf{r}_{i} \right \| /R \right ) $. The Wendland compact radial basis functions \cite{wendland1995piecewise} are shown in Table~\ref{tab.1:rbf_functions}, where, $ \eta = \left ( \left \| \textbf{r} - \textbf{r}_{i} \right \| /R \right )$. The lower-order basis functions reduce the interpolation accuracy while the higher-order basis functions require a greater computational cost. Considering this, the Wendland C2 basis function was chosen due to it providing improved smoothness in comparison to the C0 function and due to it being more efficient than the C4 and C6 functions and takes the form of  
\begin{equation}
\label{eq:2}
\phi(\eta) = \left \{ \begin{array}{ll}
\left ( 1 - \eta \right )^{4} \left ( 4 \eta + 1 \right ) & 0 \leq \eta < 1\\ 
0 & \eta \geq 1
\end{array} \right . \, .
\end{equation}

\section{Multi-Level Greedy Surface Point Selection}
\label{sec:data_reduction}
The high-quality mesh deformation properties of radial basis functions make them appealing, however, their relatively high computational cost may prevent their use on larger-scale problems where the number of surface points can frequently exceed $10^{5}$. As introduced in equations~\ref{eq:9}~\&~\ref{eq:10}, the size of the surface and volume basis matrices are $N_s^2$ and $N_s \times N_v$ respectively. The relative CPU cost associated with solving the linear algebra in equation~\ref{eq:9} thus scales with $N_s^3$ while the CPU cost of interpolating equation~\ref{eq:10} scales with $N_s \times N_v$. Methods throughout the literature have identified the size of $N_s$ as being an issue for the solving of equation~\ref{eq:9} and the interpolation of equation~\ref{eq:10}. Notably, Rendall and Allen published a method for reducing the number of surface points based on a greedy algorithm \cite{rendall2010reduced}. Their method starts with an initial control point and sequentially uses radial basis function interpolation to find the subsequent control point with the largest error signal. This control point is then added to the next step for the process to then repeat itself until the interpolation meets a sufficient tolerance. This greedy point selection will now be introduced. 

An initial control point is first selected and added to the control points vector, $\textbf{X}_{c}$. The choice of the first control point is arbitrary thus, the first point on the list of surface nodes is used, $\textbf{X}_{s}^{(0)}$. Such that

\begin{equation}
    \textbf{X}_{c} = \textbf{X}_{s}^{(0)} \, ,
\end{equation}
where the subscript ``$c$'' denotes the control points and the subscript ``$s$'' represents the surface points. The size of the control points vector is of the size $N_{c}$, and since this is the first iteration of the greedy selection process $\textbf{X}_{c}$ contains only a singular element.

An error signal is then used to guide the greedy algorithm when selecting control points. Within this work an error vector, $\textbf{E}$, based on the difference between actual surface displacements and computed surface displacements is used

\begin{equation}
\label{eq:12}
    \textbf{E} = \Delta S - \boldsymbol \Phi_{s,c} \boldsymbol \alpha \, ,
\end{equation}
with the basis matrix, $\boldsymbol \Phi_{s, \, c}$, now being of the size of $N_s \times N_c$ and where $\textbf{E}$ is of the size of $N_s$. The element of the error vector with the largest error signal, $\textbf{E}^{(\textrm{max})}$, is then used to identify the subsequent control point to be added to $\textbf{X}_{c}$. This greedy selection process continues until the largest error signal normalized by $\Delta S$ meets a desired tolerance, $\varepsilon$, as described by

\begin{equation}
\dfrac{ \textbf{E}^{(\textrm{max})}}{\Delta S} < \varepsilon \, .
\end{equation}
The size of $\textbf{X}_{c}$ therefore depends on the number of iterations of the greedy selection process. 

The greedy point selection thus allows for the reduction in the number of points being used for solving the linear algebra in equation~\ref{eq:9} and for the interpolation in equation~\ref{eq:10}, however, the greedy system itself requires significant outlay as every time a control point is added to the system the coefficients need to be solved once more. This therefore means that the CPU cost associated with solving the greedy selection process becomes of the order of $N_c^4$. Wang et al. addressed this issue by introducing a multi-level subspace radial basis function interpolation where, at the end of each level, the error of that interpolation step is used as the object for the subsequent interpolation step \cite{wang2015improved} and can be expressed as

\begin{equation}
    \Delta S_{l+1} = \textbf{E} \, ,
\end{equation}
where the subscript ``$l+1$'' denotes the next level of the multi-level greedy selection process. The new error is then computed based on the residual from the previous step and with $ \Delta S_{l+1} \ll \Delta S_{l}$ the size of the displacements is reduced significantly. The error signal described by equation~\ref{eq:12} is thus updated error and becomes

\begin{equation}
     \Delta S_{l+1} = \Delta S_{l} - \boldsymbol \Phi_{s, c} \boldsymbol \alpha \, .
\end{equation}
The number of subspace levels, $N_l$, is used to continuously reduce the error which allows the multi-level greedy method to be more efficient than the single-level greedy method. Since the CPU cost for the multi-level greedy method is of the order of $N_l \times N_c^4$ whereas the single-level greedy method is of the order of $(N_l \times N_c)^4$. The multi-level greedy selection process can be summarised as follows 

\begin{equation}
    \begin{array}{rcccl}
        \Delta S & = & \sum_{i=0}^{i=N_l - 1} \Delta S^{(i)} & = & \sum_{i=0}^{i=N_l - 1} \boldsymbol \Phi_{s, c}^{(i)} \boldsymbol \alpha^{(i)}
    \end{array} \, ,
\end{equation}

\begin{equation}
\label{eq:19}
    \begin{array}{rcccl}
        \Delta V & = & \sum_{i=0}^{i=N_l - 1} \Delta V^{(i)} & = & \sum_{i=0}^{i=N_l - 1} \boldsymbol \Phi_{v, c}^{(i)} \boldsymbol \alpha^{(i)}
    \end{array} \, .
\end{equation}

\section{Volume Point Reduction}
\label{sec:volume}
The reduced control point list after the multi-level greedy point selection now means that the CPU cost of interpolating equation~\ref{eq:10} scales with $N_l \times N_c \times N_v$. To efficiently obtain the volume point displacements, it is thus of interest to reduce $N_v$ for large-scale problems since it can often be in the order of magnitude of $N_s \times \sqrt{N_s}$. A wall distance-based function, $\psi$, is therefore introduced to restrict the $N_v$ based on the work from Xie and Liu \cite{xie2017efficient}, namely

\begin{equation}
    \psi = \psi \left ( \frac{d(\textbf{r})}{D} \right ) \, ,
\end{equation}
where $d (\textbf{r})$ signifies the wall distance and $D$ is the support distance of the wall function. The wall distance function is of compact support, which means it decays from away from the wall and is zero outwith the support distance as shown,

\begin{equation}
    \psi(\xi) = \left \{ \begin{array}{ll}
    \left ( 1 - \xi \right ) & 0 \leq \xi < 1\\ 
    0 & \xi \geq 1
    \end{array} \right . \, ,
\end{equation}
where $ \xi = d(\textbf{r}) / D $ is the wall distance normalised by the support distance. The distance $D$ is computed as a function of the maximum surface displacement and by using a volume reduction factor, $k$, with the subsequent expression

\begin{equation}
    D = k (\Delta S_l)^{max} \, . 
\end{equation}
The wall distance function is then included in the interpolation function, $f(\textbf{r})$, in equation~\ref{eq:1} and the updated function has the form:

\begin{equation}
\label{eq:23}
    f(\textbf{r}) =  \psi \left ( \frac{d(\textbf{r})}{D} \right ) \sum_{i=1}^{N} \alpha_{i} \phi \left ( \left \| \textbf{r} - \textbf{r}_{i} \right \| \right ) \, ,
\end{equation}
the volume points are then displaced after the first level, where $N_v$ remains relatively high to incorporate large deformations. That said, the number of control points $N_c$ remains relatively low keeping the cost of interpolation low. The volume point interpolation can be given as

\begin{equation}
    \Delta V_l = \boldsymbol \Phi_{v, c} \boldsymbol \alpha \, ,
\end{equation}
at the end of each level the support distance is updated so the interpolation target is set to the error of the previous step as shown

\begin{equation}
    D = k (\Delta S_{l+1})^{max} \, .
\end{equation}
Updating the support distance after each level is possible since $\Delta S_{l+1} \ll \Delta S_l$, meaning the number of volume points in the sphere of influence of the basis function is less. Therefore $D_{l+1} \ll D_l$ and so the number of updated volume points $N_{v,l+1} \ll N_{v,l}$ helping to reduce the cost of the interpolation step.

\section{Results}
\label{sec:results}
In order to assess the efficiency of radial basis function mesh deformation techniques for aircraft icing two separate test cases are setup. The first test case is of a two-dimensional airfoil and the second test case is of a three-dimensional swept-wing. The performance of the data-reduction schemes are then evaluated through comparisons with well-known benchmark test cases. All simulations are run on a single processor to focus on the efficiency of the mesh deformation techniques opposed to the scalability and parallelism of the code. The processor used was an Intel(R) Xeon(R) X5660 CPU with a clockspeed of 2.80GHz. 

\subsection{2D -- NACA0012 Airfoil}
\label{sec:airfoil}
The aim of the first test case is to assess the efficiency and robustness of radial basis function mesh deformation data-reduction schemes on a two-dimensional icing problem. Experimental icing tests on a NACA0012 airfoil performed in the NASA Lewis Icing Research Tunnel (IRT) \cite{ruff1990users} are used as a reference for the numerical icing predictions. The tests primarily concentrated on the repeatability of ice shapes over a range of icing conditions. One subset of these conditions is chosen and is outlined in Table~\ref{tab:2}. The test condition extends for two-minutes of ice accretion during the glaze ice regime. Ice is accreted over a NACA0012 airfoil with a chord of $0.3$ m at an angle-of-attack of $0^{\circ}$. The freestream velocity is $129$ m/s. The outside air temperature is $260.5$ K and the static pressure is $90,700$ Pa. The mean volume diameter of the supercooled water droplets is $20$ $\mu \textrm{m}$ which contain a liquid water content of $0.5$ $\textrm{g}/\textrm{m}^{3}$.   

\begin{figure}[htb!]
\centering
\subfloat[Structured volume mesh.\label{fig.1:subfig-1:airfoil_mesh}]
{
\includegraphics[height=0.31\linewidth, trim={0.5cm 6cm 3cm 4cm},clip]{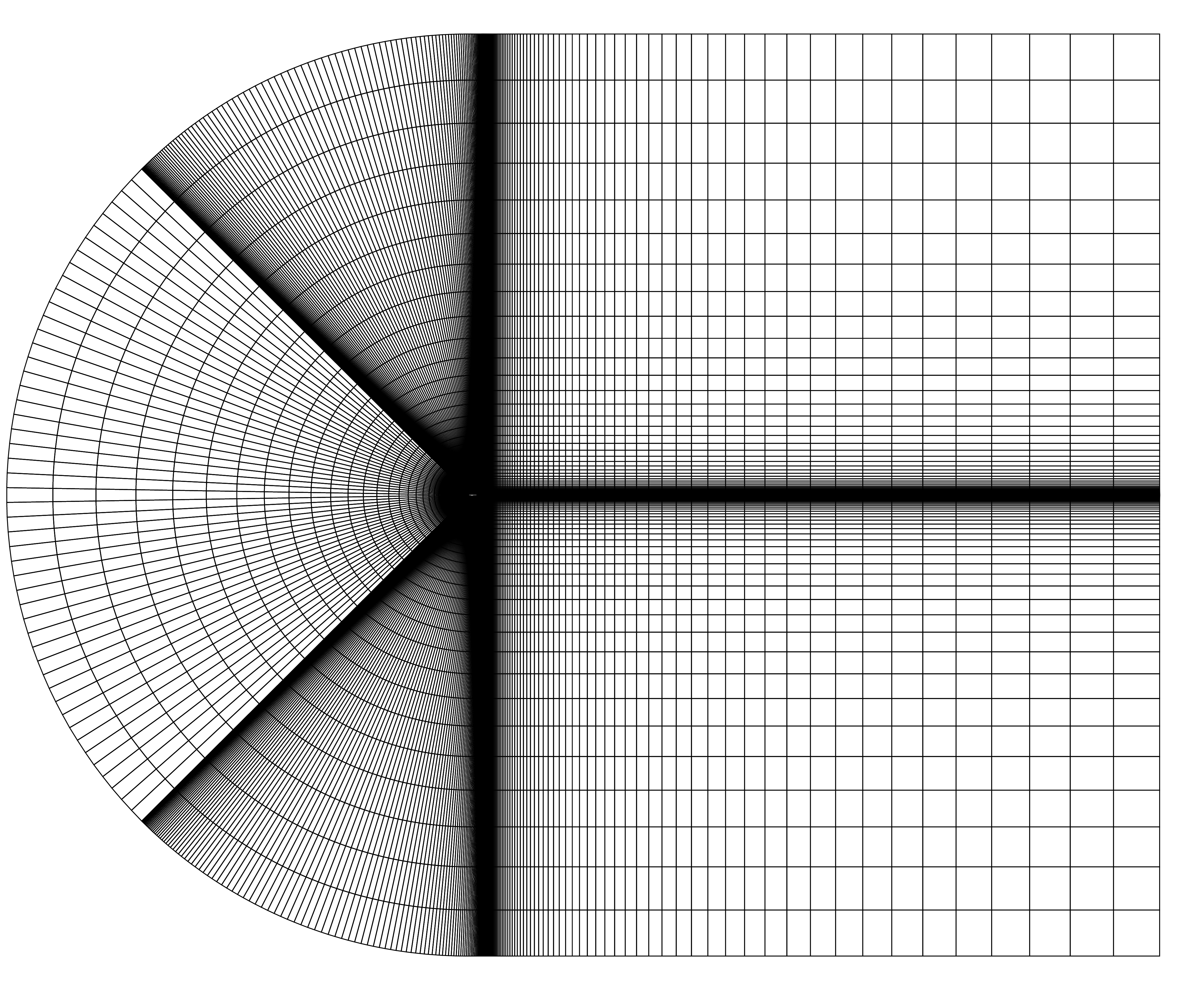}
} \hfill
\subfloat[Close-up view of the mesh surrounding the airfoil.\label{fig.1:subfig-2:airfoil_mesh}]
{
\includegraphics[height=0.31\linewidth]{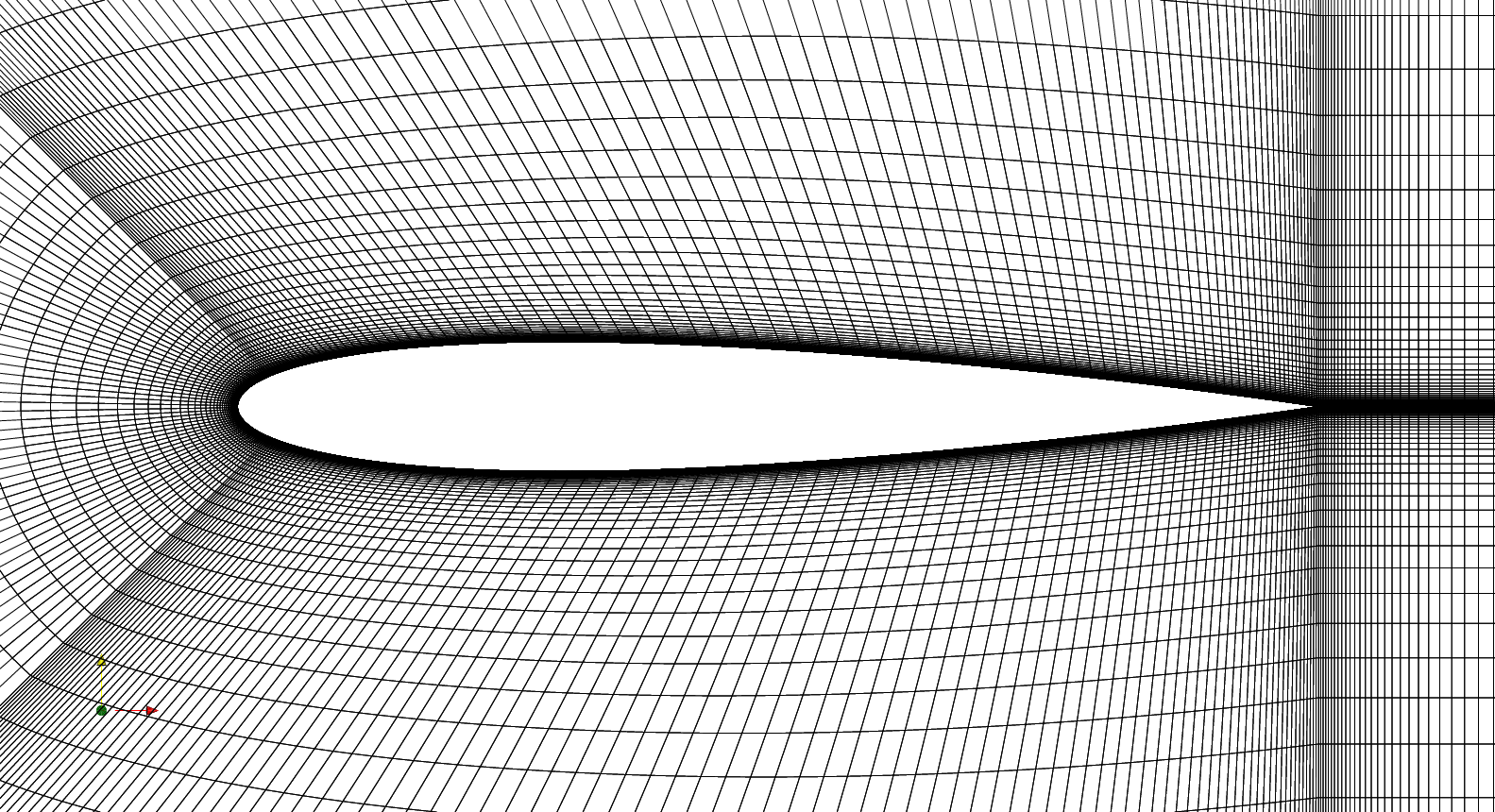}
}
\caption[Structured NACA0012 airfoil mesh]{Structured NACA0012 airfoil mesh. Constructed using the multi-block grid strategy.}
\label{fig.1:airfoil_mesh}
\end{figure}

The spatial discretization of the two-dimensional airfoil mesh is achieved using a structured multi-block grid and is shown in Fig.~\ref{fig.1:airfoil_mesh}. The total number of elements in the domain is $44,055$ and the total number of vertices is $44,500$. The NACA0012 airfoil has 247 vertices distributed around its surface with vertices congregated around the leading and trailing edges. The far-field is placed 25 chord lengths from the airfoil. The grid resolution is sufficient at the wall to ensure $\textrm{y}+ < 1$.

The standard one-equation SA turbulence model is used with the RANS equations in SU2. The SA turbulence variable is convected using a first-order scalar upwind method. The convective fluxes are computed using the Roe scheme and second-order accuracy is achieved using the MUSCL scheme. The viscous fluxes are approximated using the weighted-least-squares numerical method. The flow solution is considered converged when there is a reduction of 6 orders of magnitude on the density residual. 

\begin{figure}[htb!]
\centering
\includegraphics[width=0.7\linewidth]{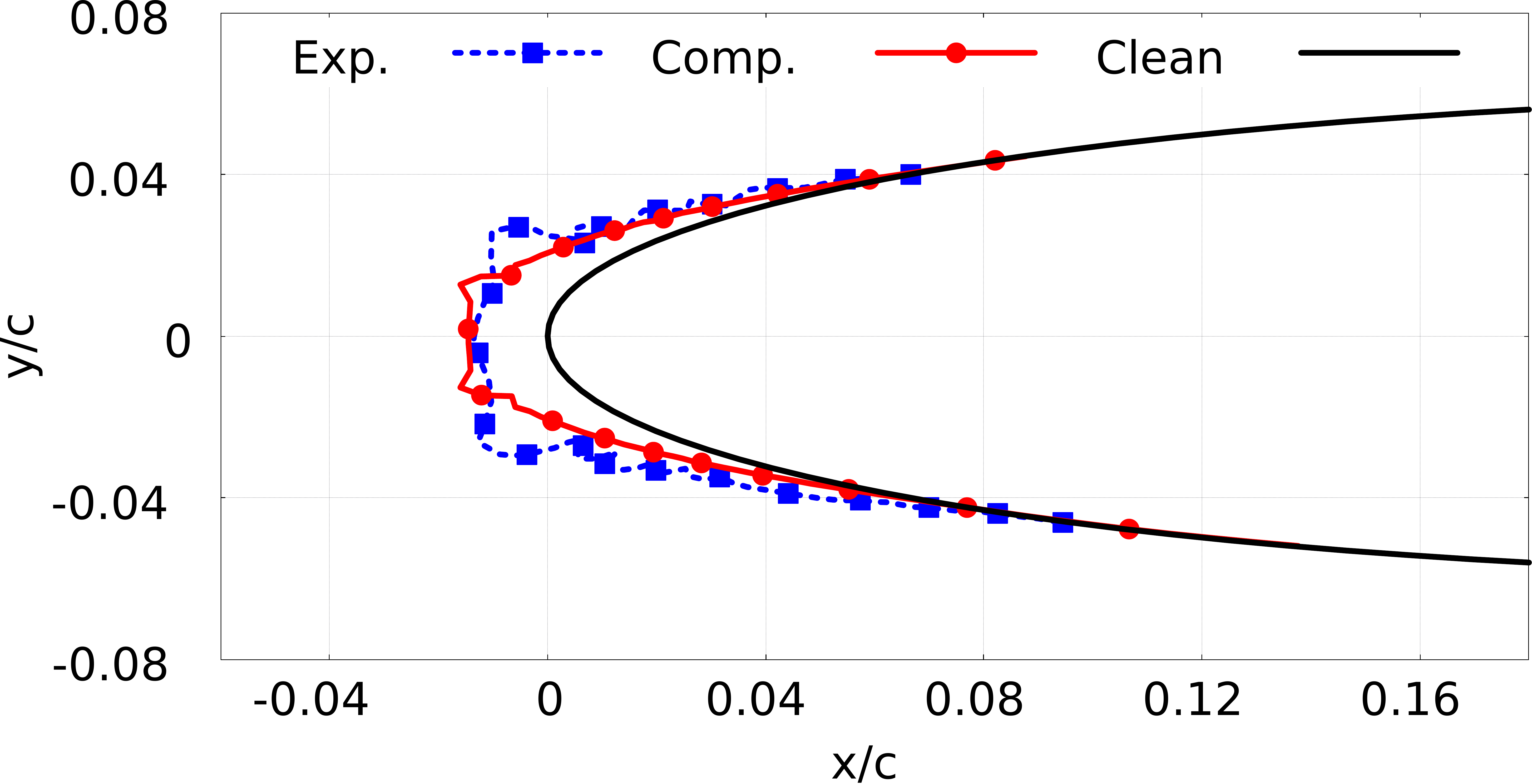}
\caption[Comparison of ice shapes on an airfoil]{Comparison of computed and experimental ice shapes on a NACA0012 airfoil under conditions identified in Table~\ref{tab:2}.}
\label{fig.2:airfoil_ice}
\end{figure}

The PoliMIce software library provides state-of-the-art ice formation models \cite{gori2015polimice}. The model used in this work to capture the complex experimental ice shapes is the local exact solution of the unsteady Stefan problem for the temperature profiles within the ice layer in glaze conditions \cite{gori2018local}. Multi-step ice accretion simulations are performed at $5$ second intervals to iteratively update the solution and account for unsteady ice accretion. The final predicted ice shape after 120 seconds of ice accretion is compared against the measured data from Ref.~\cite{ruff1990users} and is shown in Fig.~\ref{fig.2:airfoil_ice}. The ice shape exhibits distinct horns paradigmatic of the glaze ice regime and the overall mass of ice is predicted well. 

To deform the iced profile the radial basis function mesh deformation uses the compactly supported Wendland C2 function. Points within the support radius of $R = 2c$ are deformed. Altogether five levels of multi-level greedy surface point reduction are used. Each level is updated when there is a reduction of $\varepsilon = 10^{-1}$ in the normalized error. A volume reduction factor of $k = 5$ is chosen. 

To assess the performance of the deformation techniques, a benchmark test case taken from Ref.~\cite{xie2017efficient} is used as a reference. The benchmark test case imposes a sinusoidal motion to the airfoil which can be described by:

\begin{equation}
\Delta y = 0.01 \sin \left ( 15 \pi x \right ) \, ,
\end{equation}
where $\Delta y$ describes the displacement of the airfoil as a function of its position in the $x$-direction along its chord, $c$. The sinusoidal motion thus represents a more globalised deformation of the airfoil. The coefficients describing the amplitude and wave length of the sinusoidal function are modified to better represent the current geometry. Identical radial basis function parameters and mesh are used to ensure maximum similarity and to allow for direct comparisons between local and global deformations.  

\begin{figure}[htb!]
\centering
\subfloat[Convergence history in terms of selected points.\label{fig.3:subfig-1:airfoil_greedyConvergence}]
{
\includegraphics[width=0.6\linewidth]{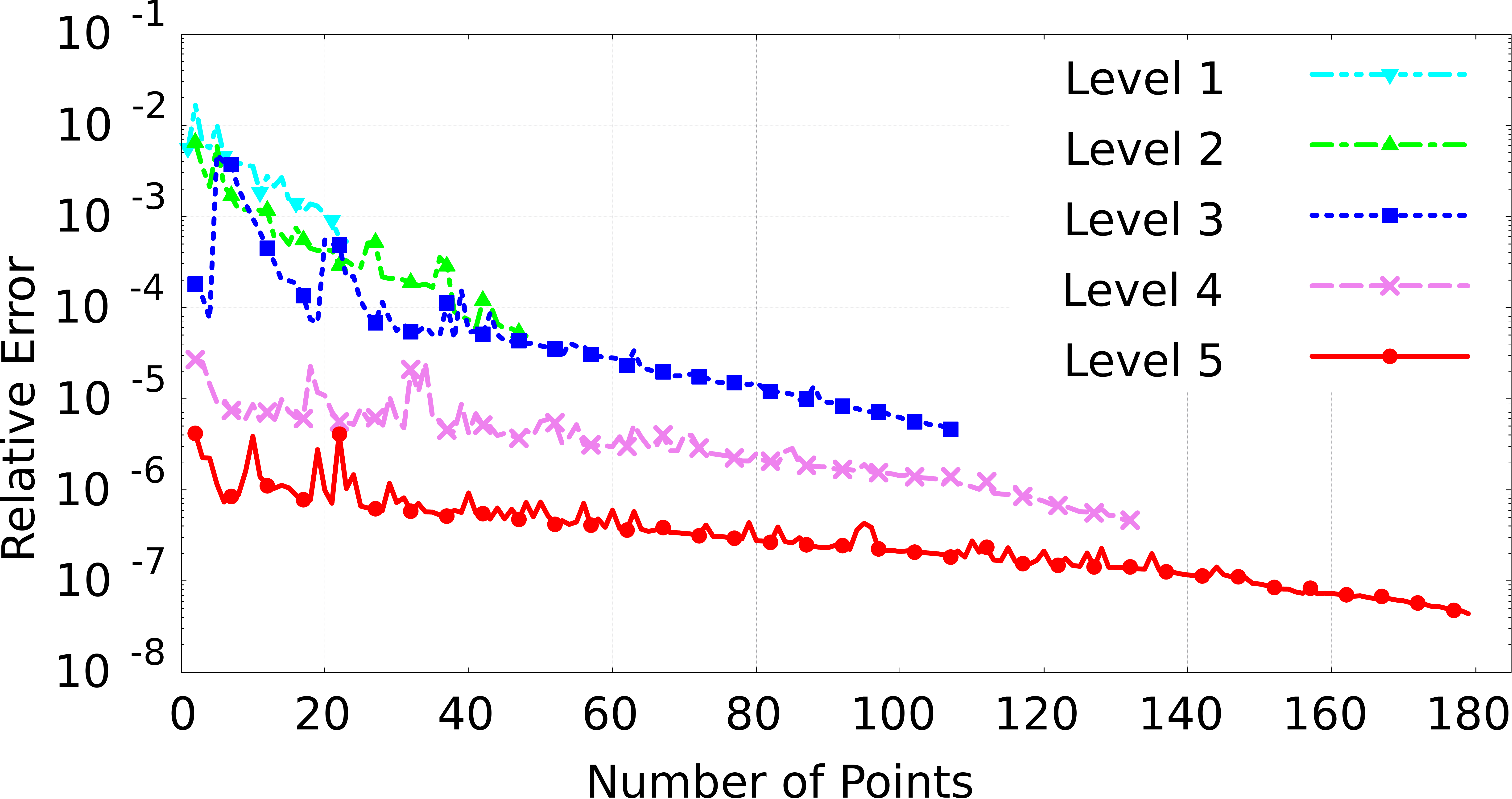}
} \\
\subfloat[Convergence history in terms of CPU time.\label{fig.3:subfig-2:airfoil_greedyConvergence}]
{
\includegraphics[width=0.6\linewidth]{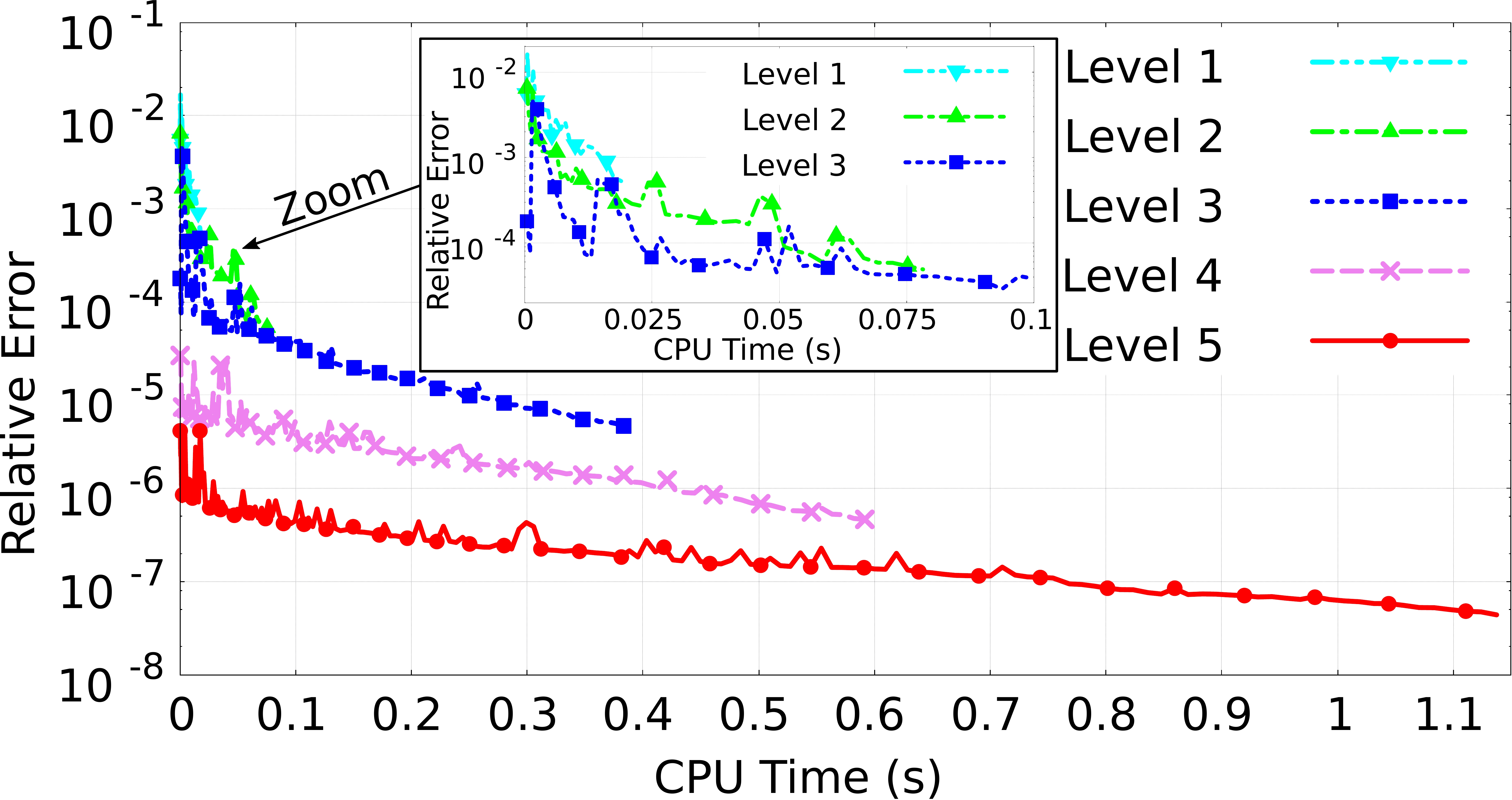}
}
\caption[Comparison of error reduction rates on the iced airfoil]{Comparison of error reduction rates in terms of selected points and CPU time for the NACA0012 airfoil.}
\label{fig.3:airfoil_greedyConvergence}
\end{figure}

\begin{figure}[htb!]
\centering
\subfloat[Convergence history in terms of selected points.\label{fig.4:subfig-1:defAirfoil_greedyConvergence}]
{
\includegraphics[width=0.6\linewidth]{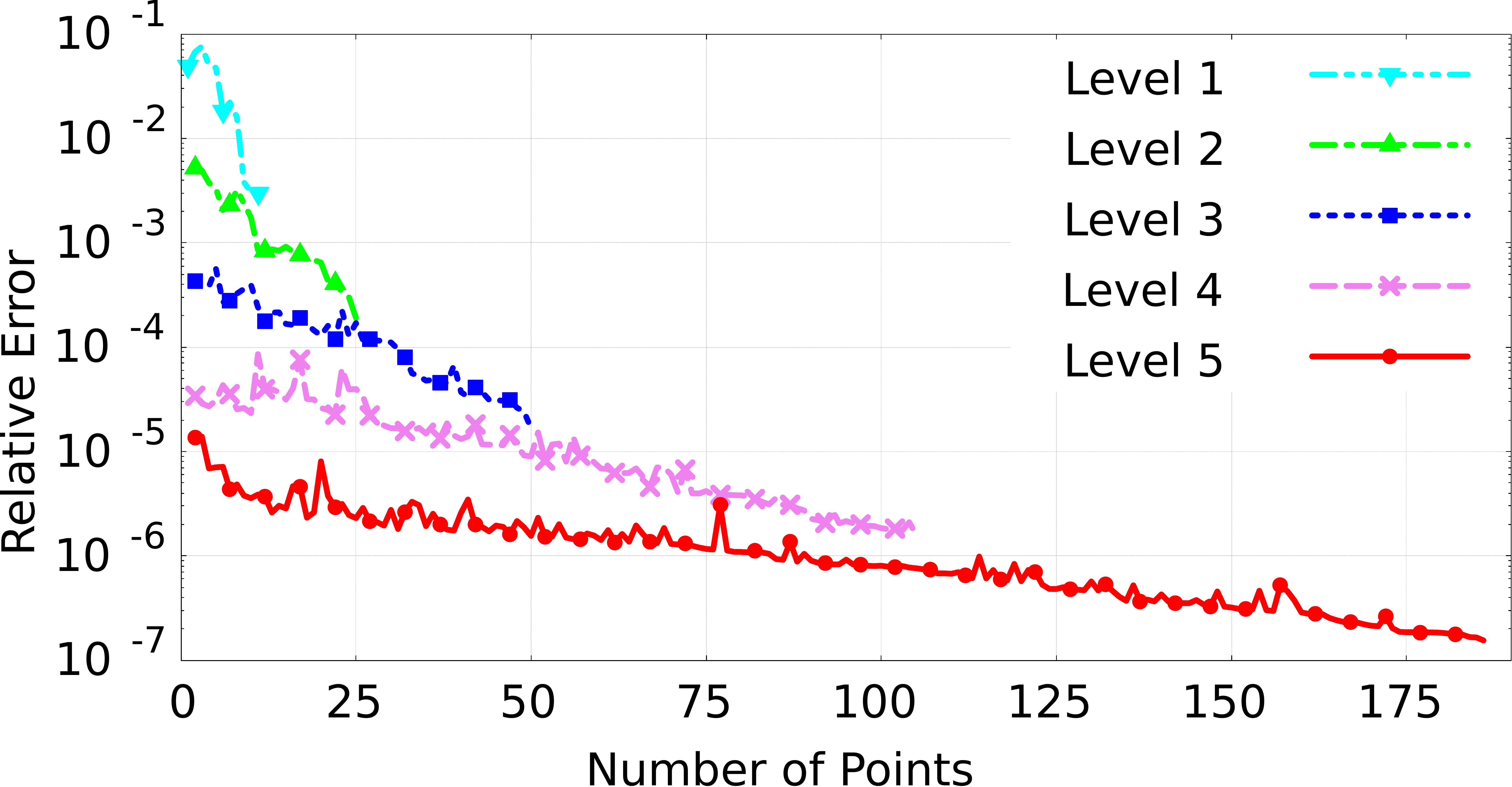}
} \\
\subfloat[Convergence history in terms of CPU time.\label{fig.4:subfig-2:defAirfoil_greedyConvergence}]
{
\includegraphics[width=0.6\linewidth]{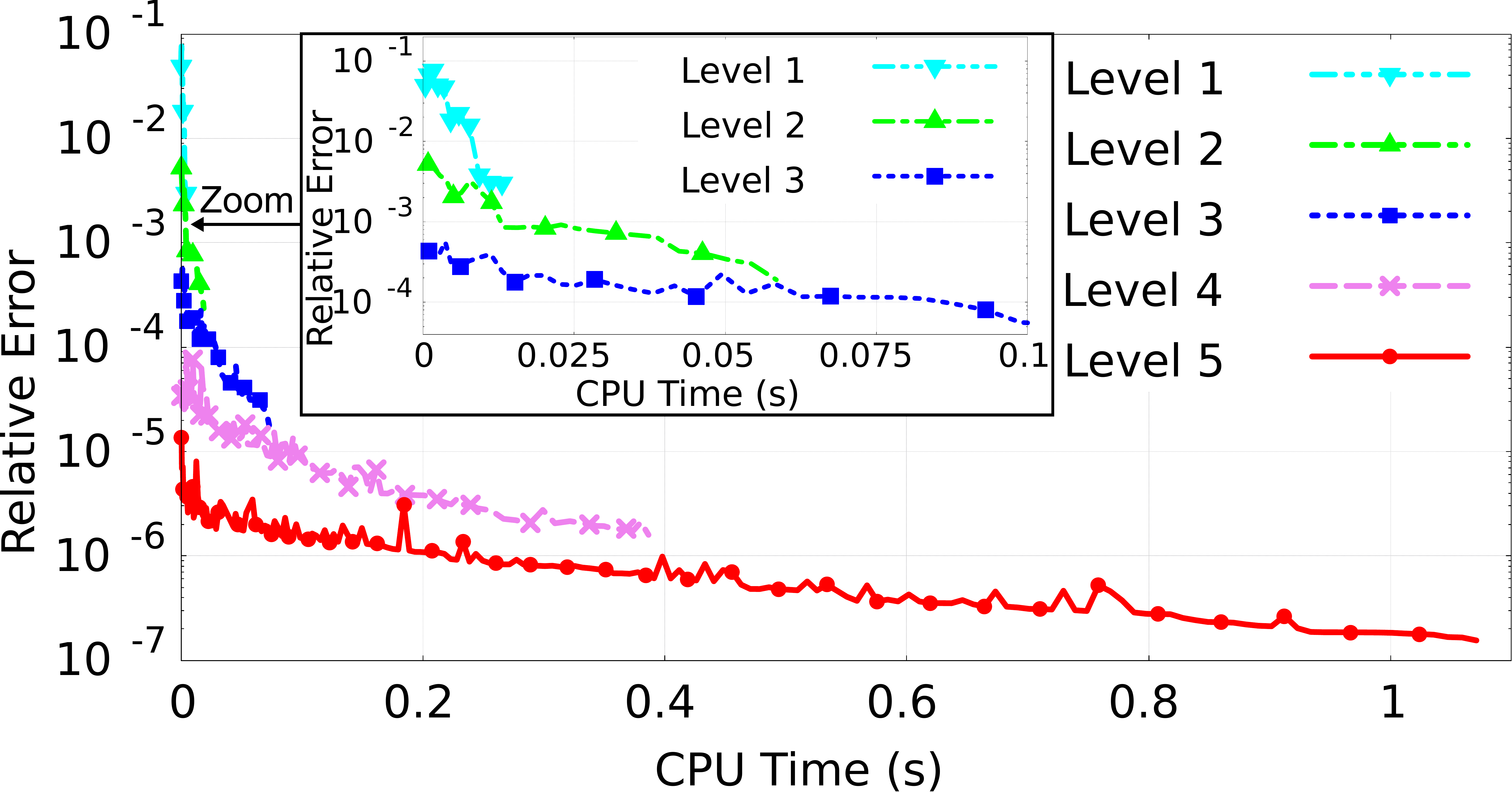}
}
\caption[Comparison of error reduction rates on the sinusiodal airfoil]{Comparison of error reduction rates in terms of selected points and CPU time for the NACA0012 airfoil with sinusoidal motion.}
\label{fig.4:defAirfoil_greedyConvergence}
\end{figure}

The convergence history of the normalized displacement error during the iced and sinusoidal deformation are respectively illustrated in Figs.~\ref{fig.3:airfoil_greedyConvergence} \&~\ref{fig.4:defAirfoil_greedyConvergence}. The efficiency of the multi-level greedy point selection performs very well in terms of maintaining a low number of control points and in terms of CPU time for both test cases. Both kinds of deformation achieve 5 levels greedy surface point selection. The sinusoidal deformation shows marginally improved performance during the initial levels. However, at the highest level both kinds of deformation show similar performance. The iced deformation shown in Fig.~\ref{fig.3:airfoil_greedyConvergence} requires 179 control points at the $5^{\textrm{th}}$ level and has an associated CPU time of approximately 1.14 s. While the sinusoidal deformation shown in Fig.~\ref{fig.4:defAirfoil_greedyConvergence} requires 182 control points at the $5^{\textrm{th}}$ level and thus exhibits a similar CPU time of approximately 1.16 s.

\begin{figure}[htb!]
\centering
\subfloat[Level 1: 23 Control Points.\label{fig.5:subfig-1:airfoil_greedySelection}]
{
\includegraphics[width=0.48\linewidth]{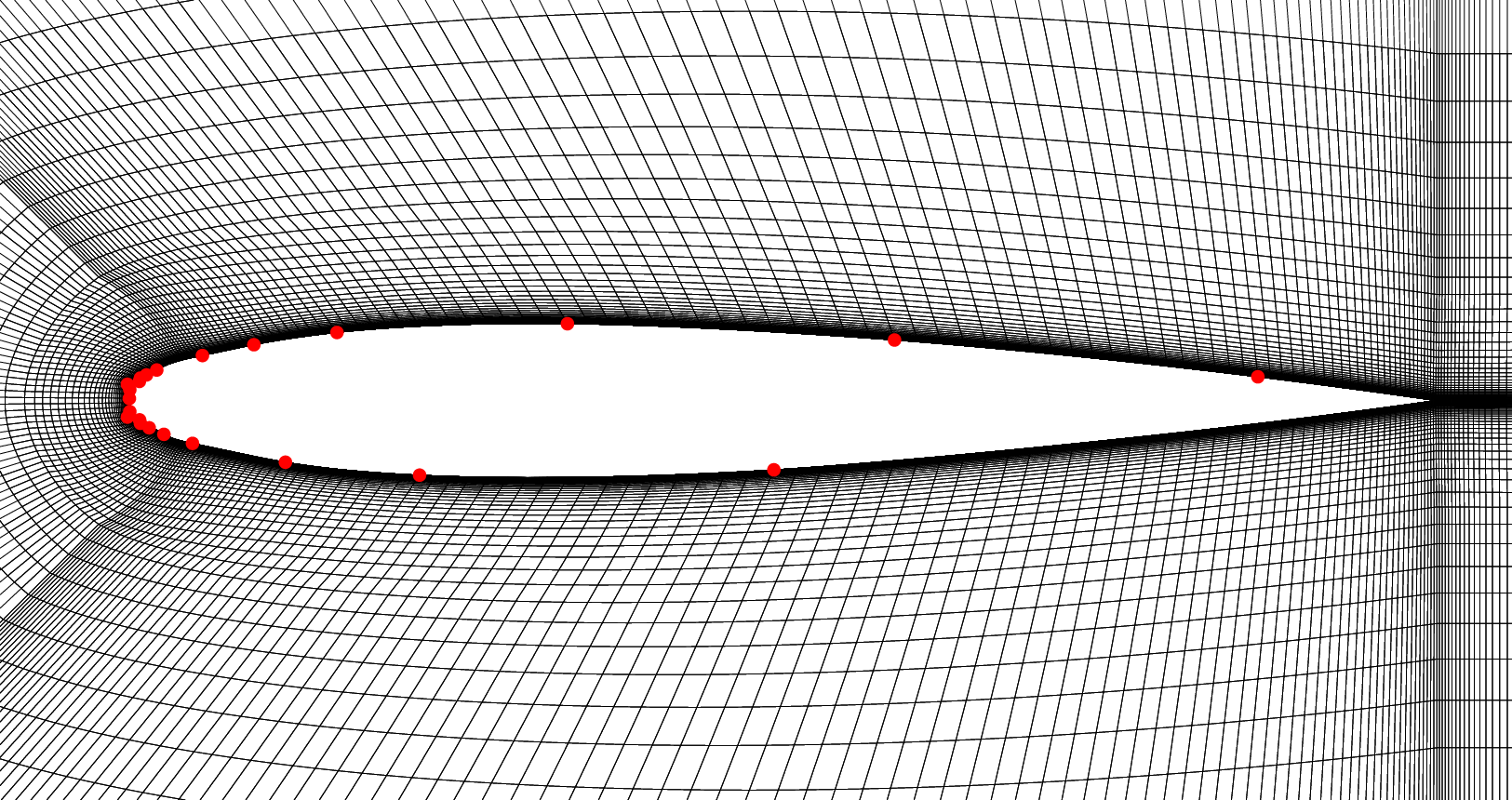}
} \hfill
\subfloat[Level 2: 48 Control Points.\label{fig.5:subfig-2:airfoil_greedySelection}]
{
\includegraphics[width=0.48\linewidth]{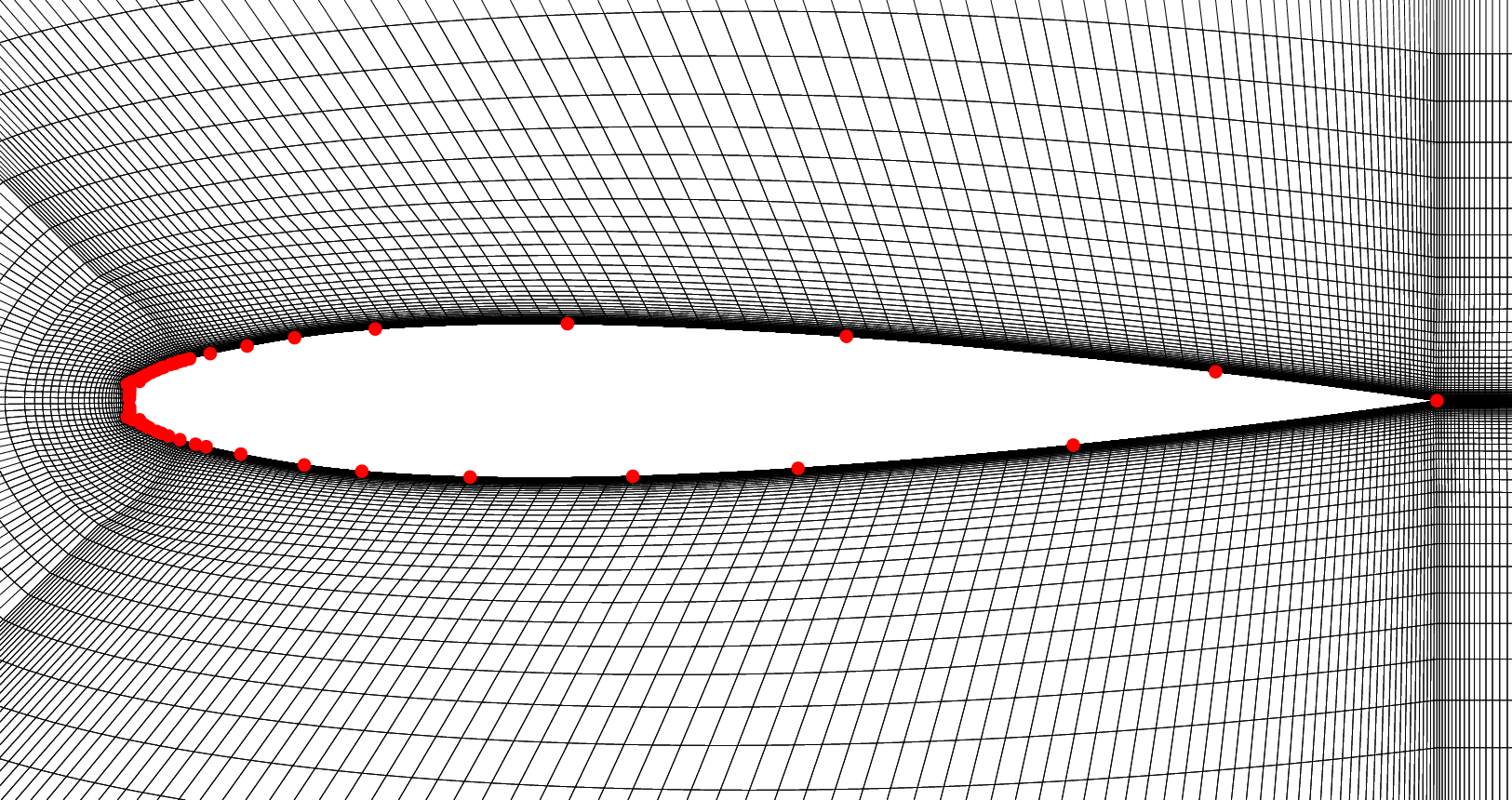}
} \\
\subfloat[Level 3: 107 Control Points.\label{fig.5:subfig-3:airfoil_greedySelection}]
{
\includegraphics[width=0.48\linewidth]{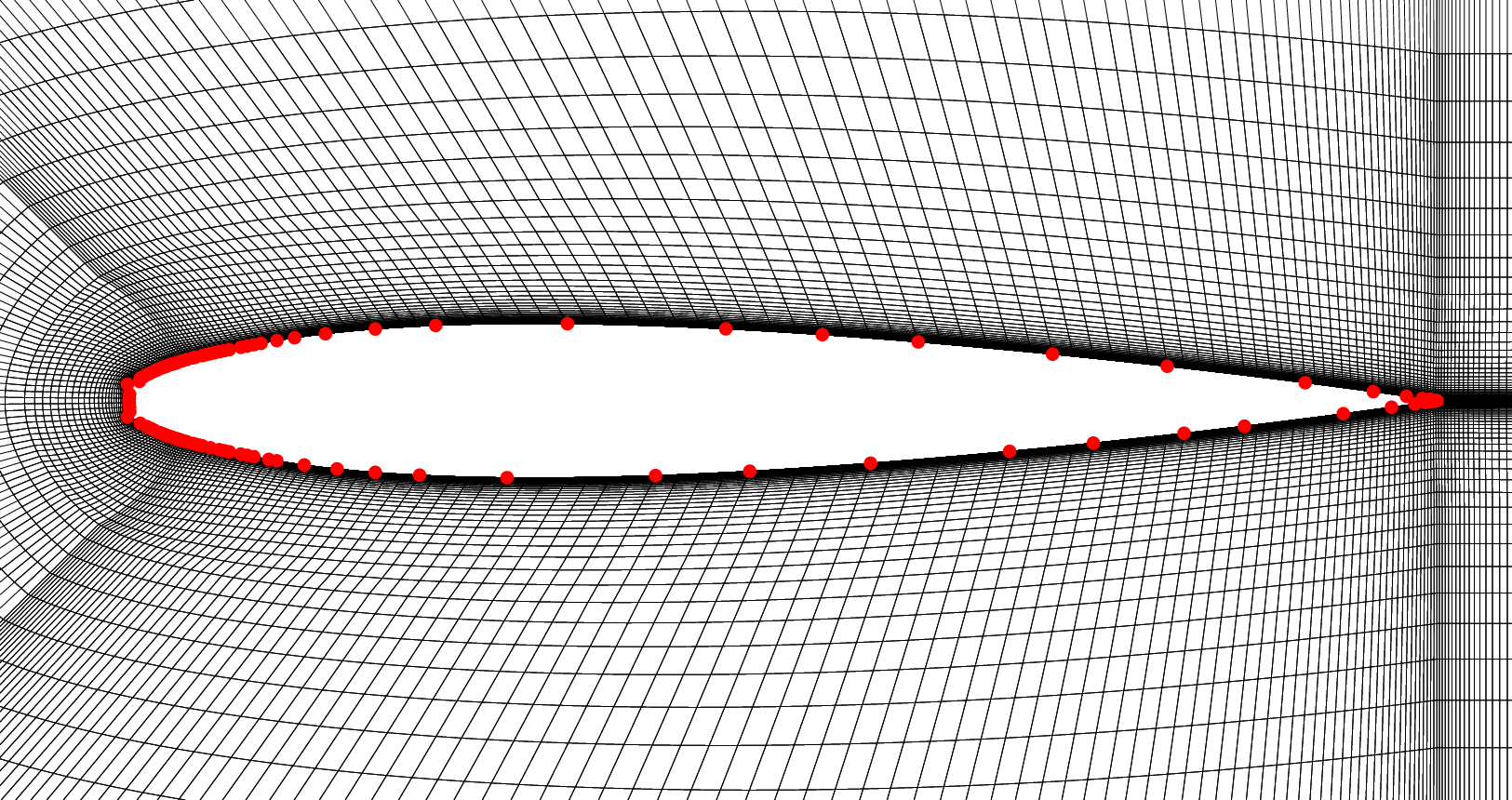}
} \hfill
\subfloat[Level 4: 132 Control Points.\label{fig.5:subfig-4:airfoil_greedySelection}]
{
\includegraphics[width=0.48\linewidth]{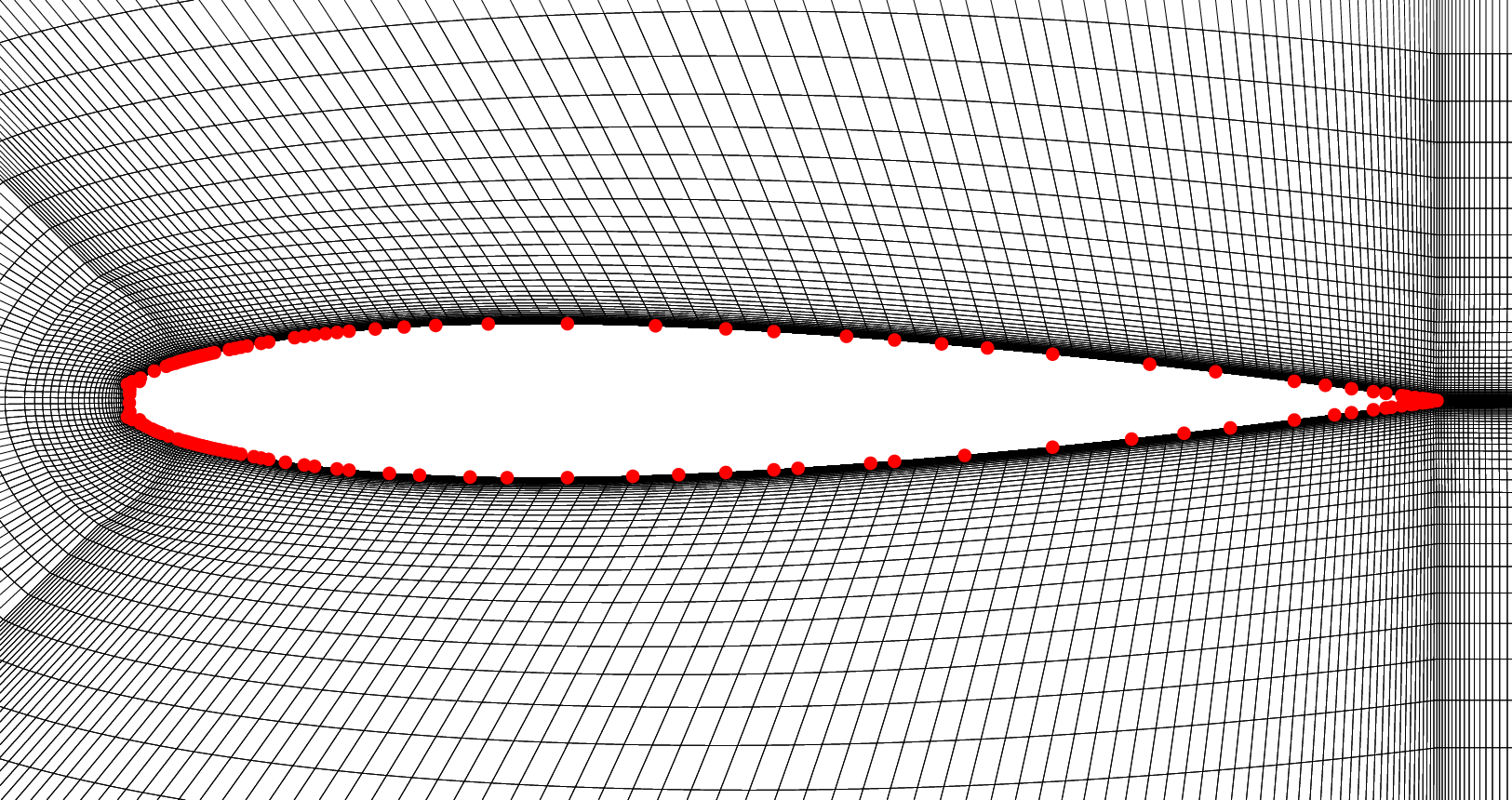}
}
\caption[Control points selected on the iced airfoil]{Control points selected during the multi-level greedy point selection for the NACA0012 airfoil; where the red points indicate the control points.}
\label{fig.5:airfoil_greedySelection}
\end{figure}

\begin{figure}[htb!]
\centering
\subfloat[Level 1: 11 Control Points.\label{fig.6:subfig-1:defAirfoil_greedySelection}]
{
\includegraphics[width=0.48\linewidth]{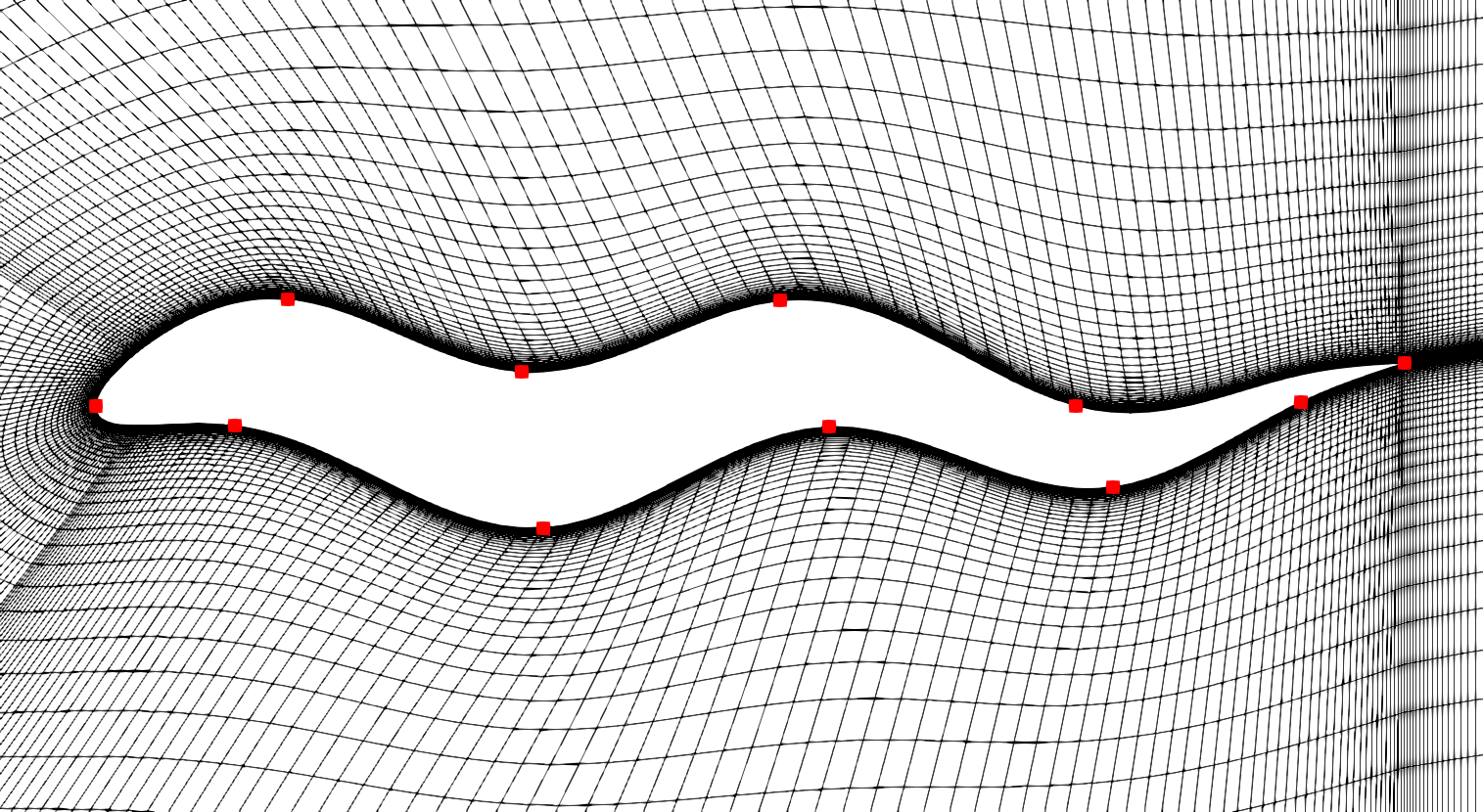}
} \hfill
\subfloat[Level 2: 27 Control Points.\label{fig.6:subfig-2:defAirfoil_greedySelection}]
{
\includegraphics[width=0.48\linewidth]{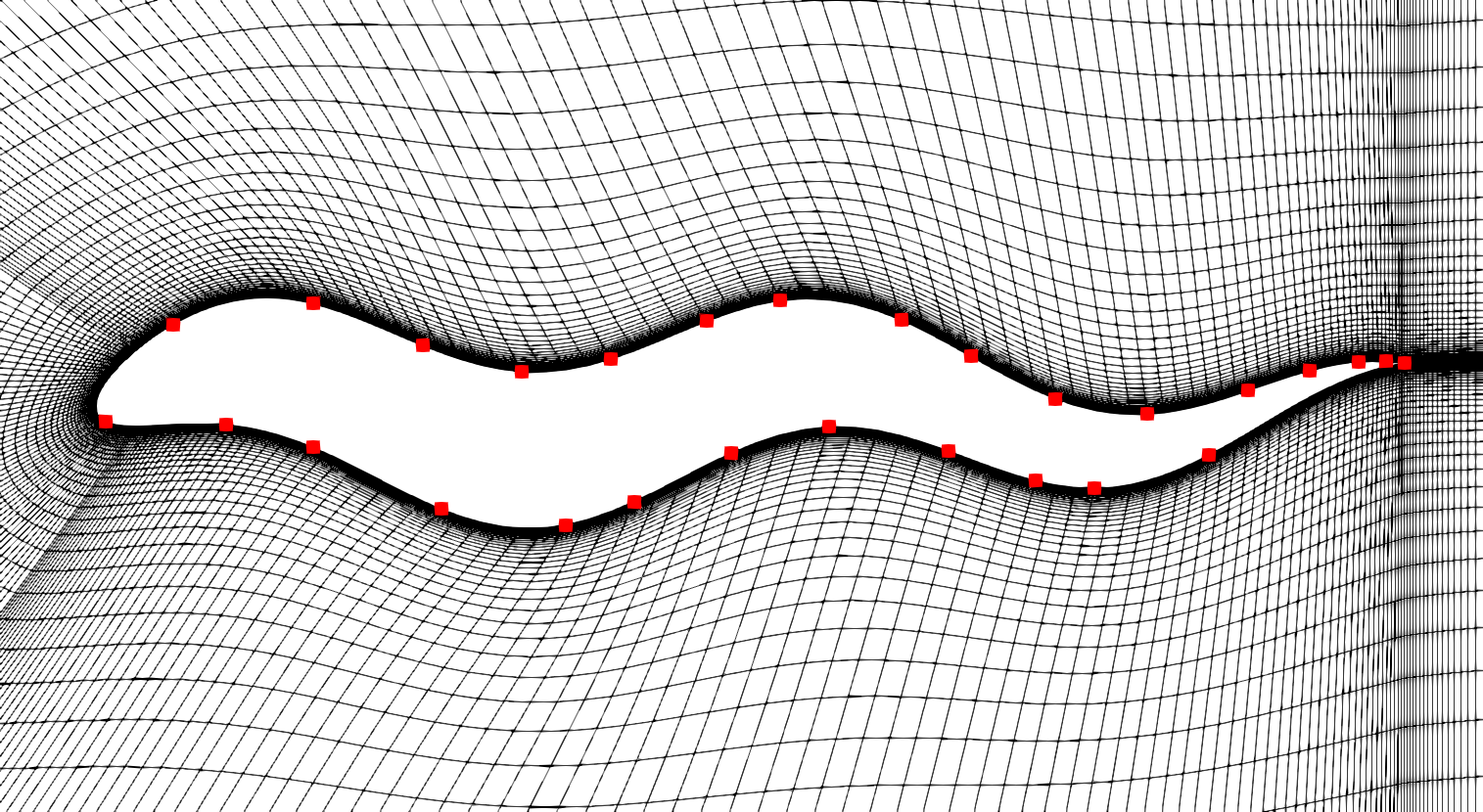}
} \\
\subfloat[Level 3: 53 Control Points.\label{fig.6:subfig-3:defAirfoil_greedySelection}]
{
\includegraphics[width=0.48\linewidth]{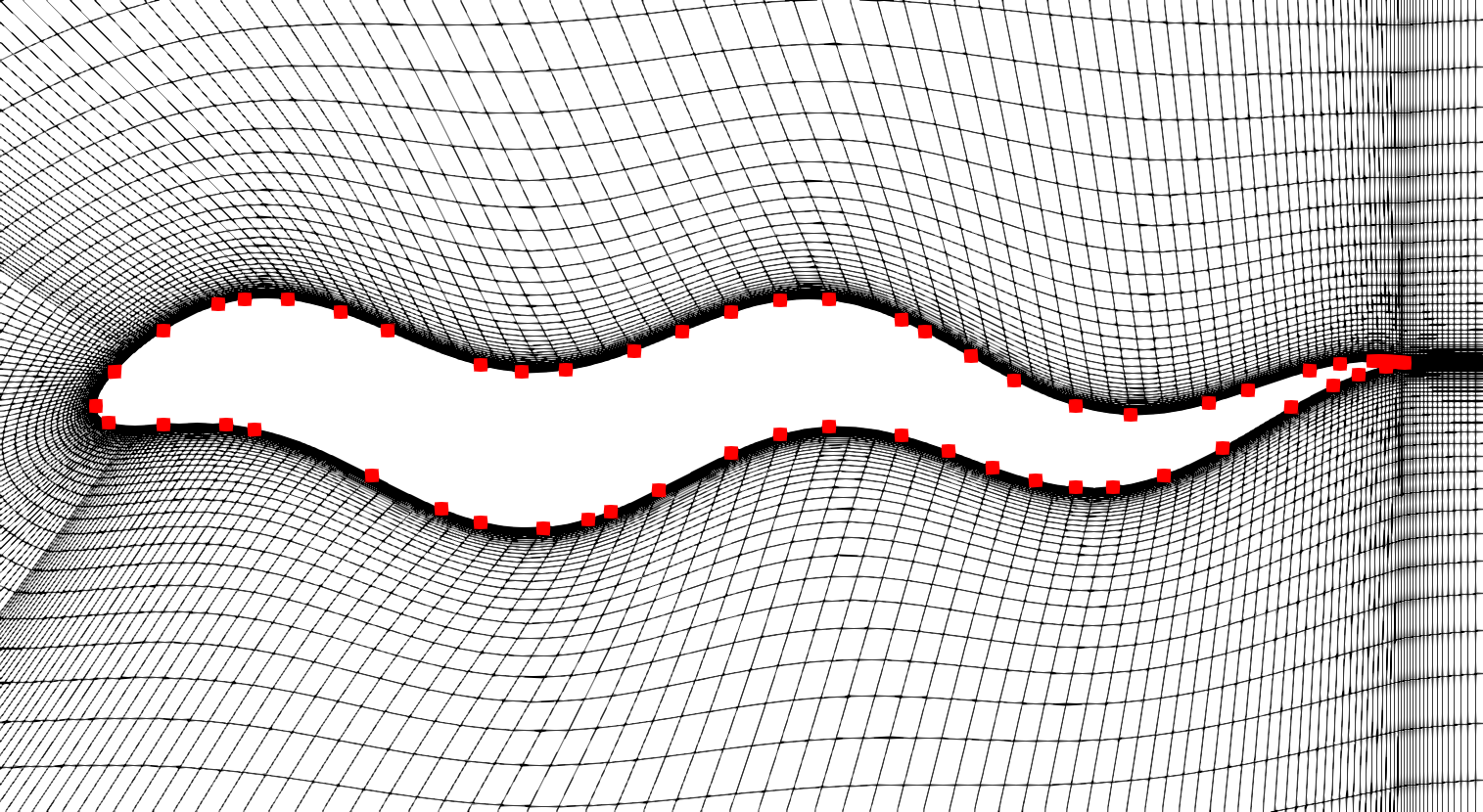}
} \hfill
\subfloat[Level 4: 108 Control Points.\label{fig.6:subfig-4:defAirfoil_greedySelection}]
{
\includegraphics[width=0.48\linewidth]{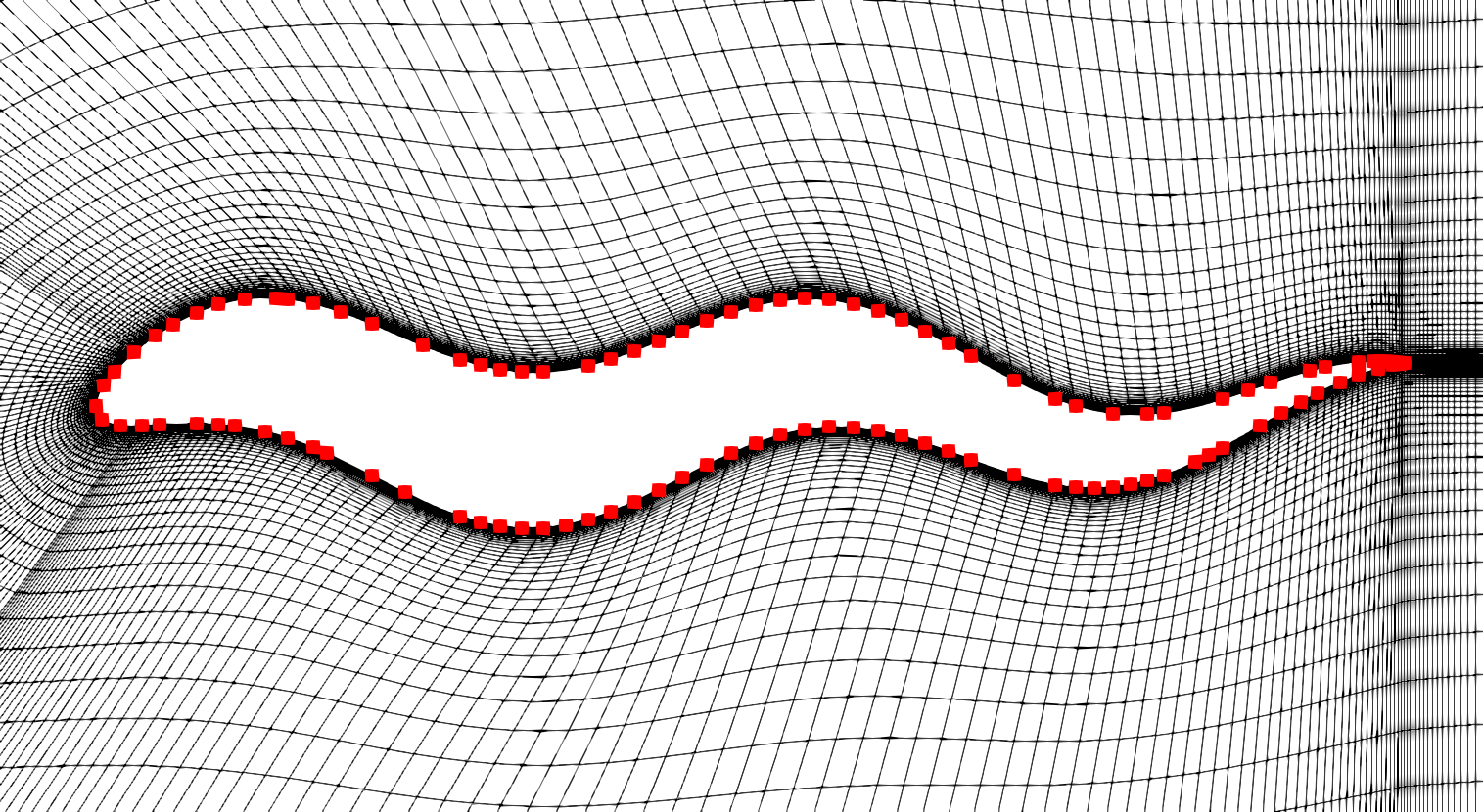}
}
\caption[Control points selected on the sinusiodal airfoil]{Control points selected during the multi-level greedy point selection for the airfoil with sinusoidal motion; where the red points indicate the control points.}
\label{fig.6:defAirfoil_greedySelection}
\end{figure}

The control points selected from the $1^{\textrm{st}}$ $\rightarrow$ $4^{\textrm{th}}$ levels of the greedy process during the iced and sinusoidal deformation are respectively shown in Figs.~\ref{fig.5:airfoil_greedySelection} \&~\ref{fig.6:defAirfoil_greedySelection}. The local ice deformation is represented by the control points which are predominately distributed around the leading edge of the airfoil as shown in Fig.~\ref{fig.5:airfoil_greedySelection}. This leads to a highly anisotropic control point distribution. It is evident that for ice shapes with any kind of level of roughness or horns a large number of control points are required otherwise the radial interpolation may potentially smooth over these local features. The global sinusoidal deformation of the airfoil is represented by the control points which are relatively evenly distributed as shown in Fig.~\ref{fig.6:defAirfoil_greedySelection}. At the lower levels, the control points are selected at the peaks of the waves. As the greedy algorithm progresses the control point list becomes more populated and control points are present all along the airfoil.

\begin{figure}[htb!]
\centering
\subfloat[Orthogonality of the clean airfoil mesh prior to deformation.\label{fig.7:subfig-1:airfoil_meshQuality}]
{
\includegraphics[width=0.48\linewidth]{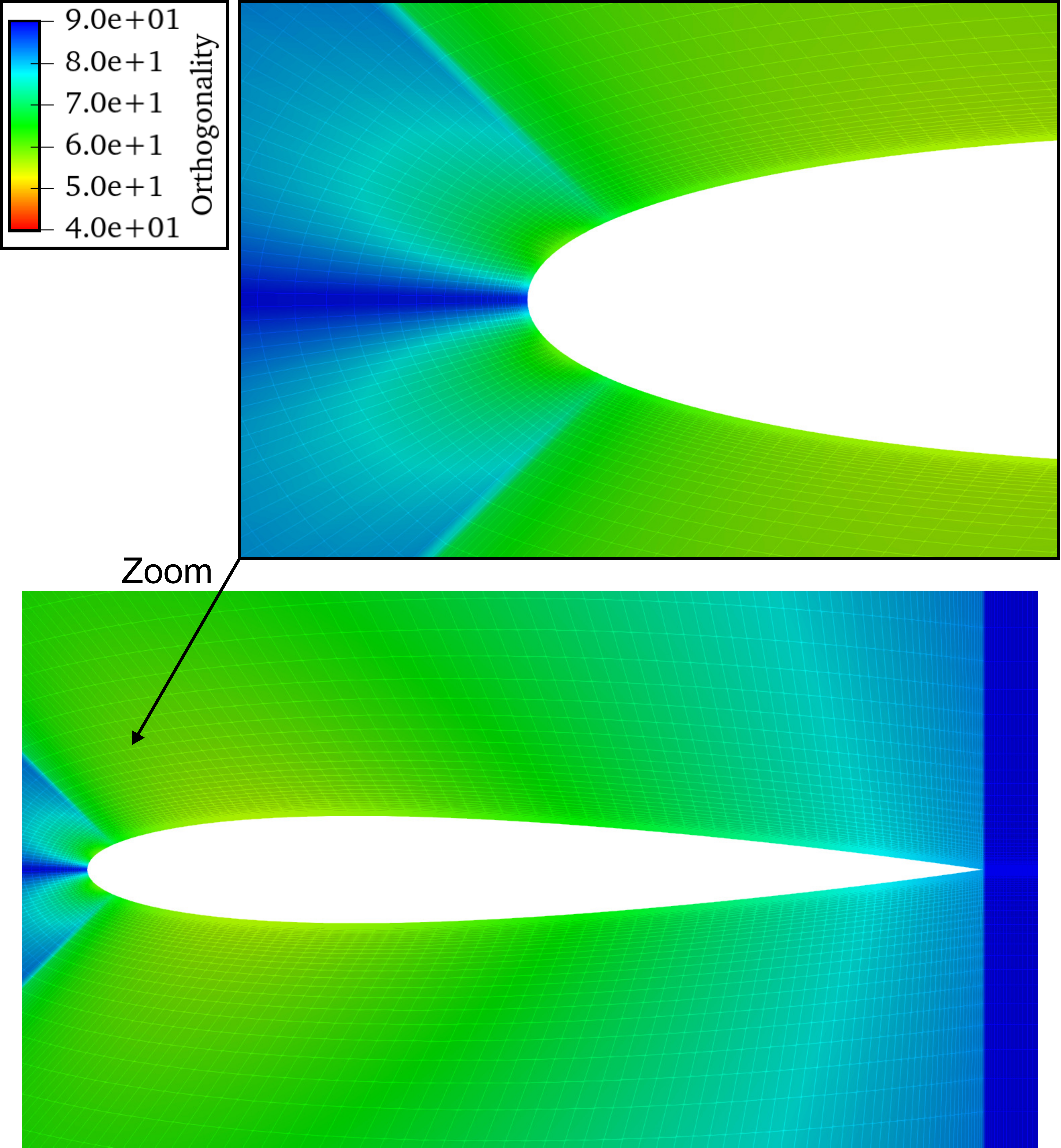}
} \hfill
\subfloat[Orthogonality of the iced airfoil mesh post deformation.\label{fig.7:subfig-2:airfoil_meshQuality}]
{
\includegraphics[width=0.48\linewidth]{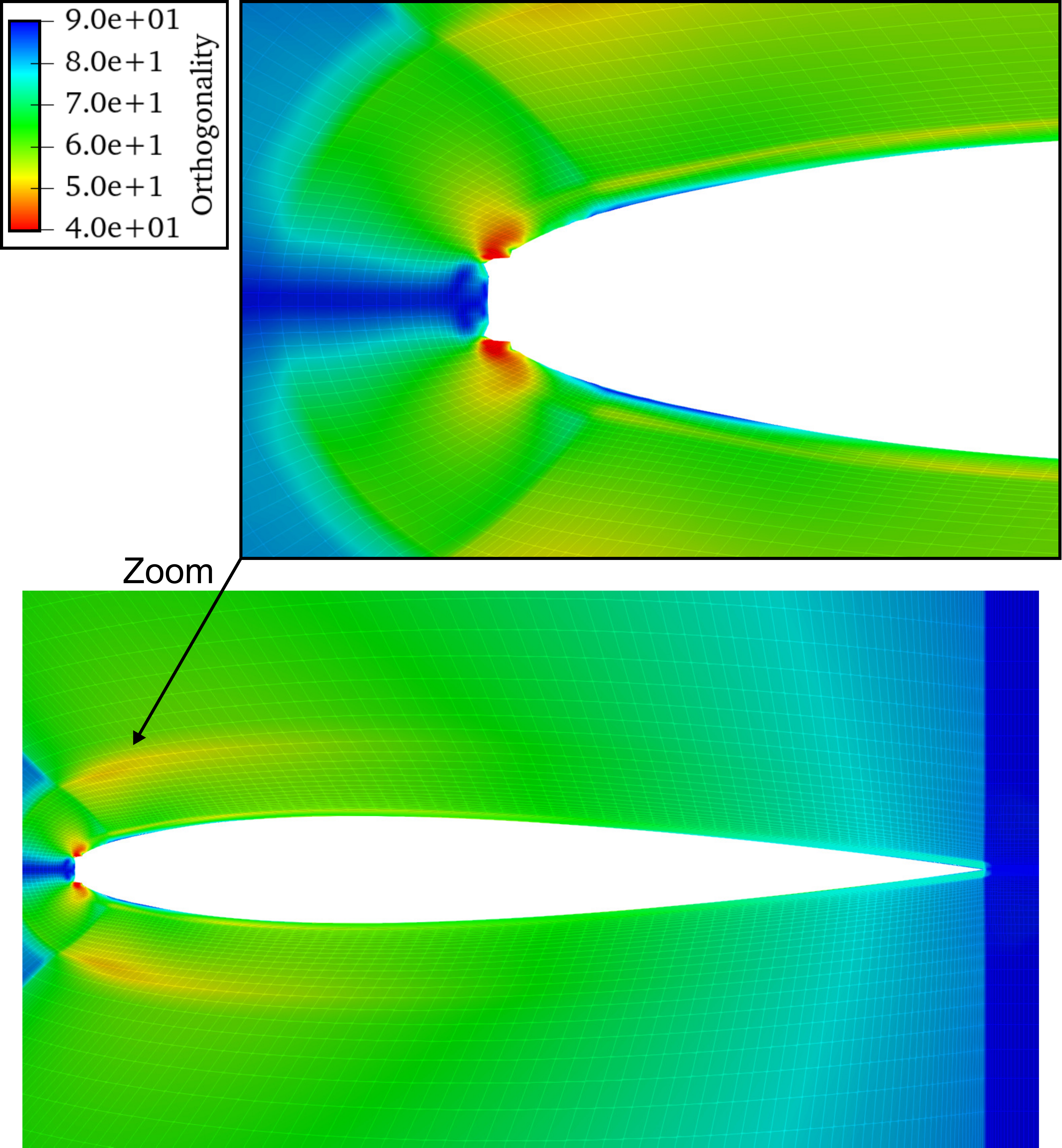}
} \\
\subfloat[Orthogonality of the sinusoidal airfoil mesh post deformation.\label{fig.7:subfig-3:airfoil_meshQuality}]
{
\includegraphics[width=0.48\linewidth]{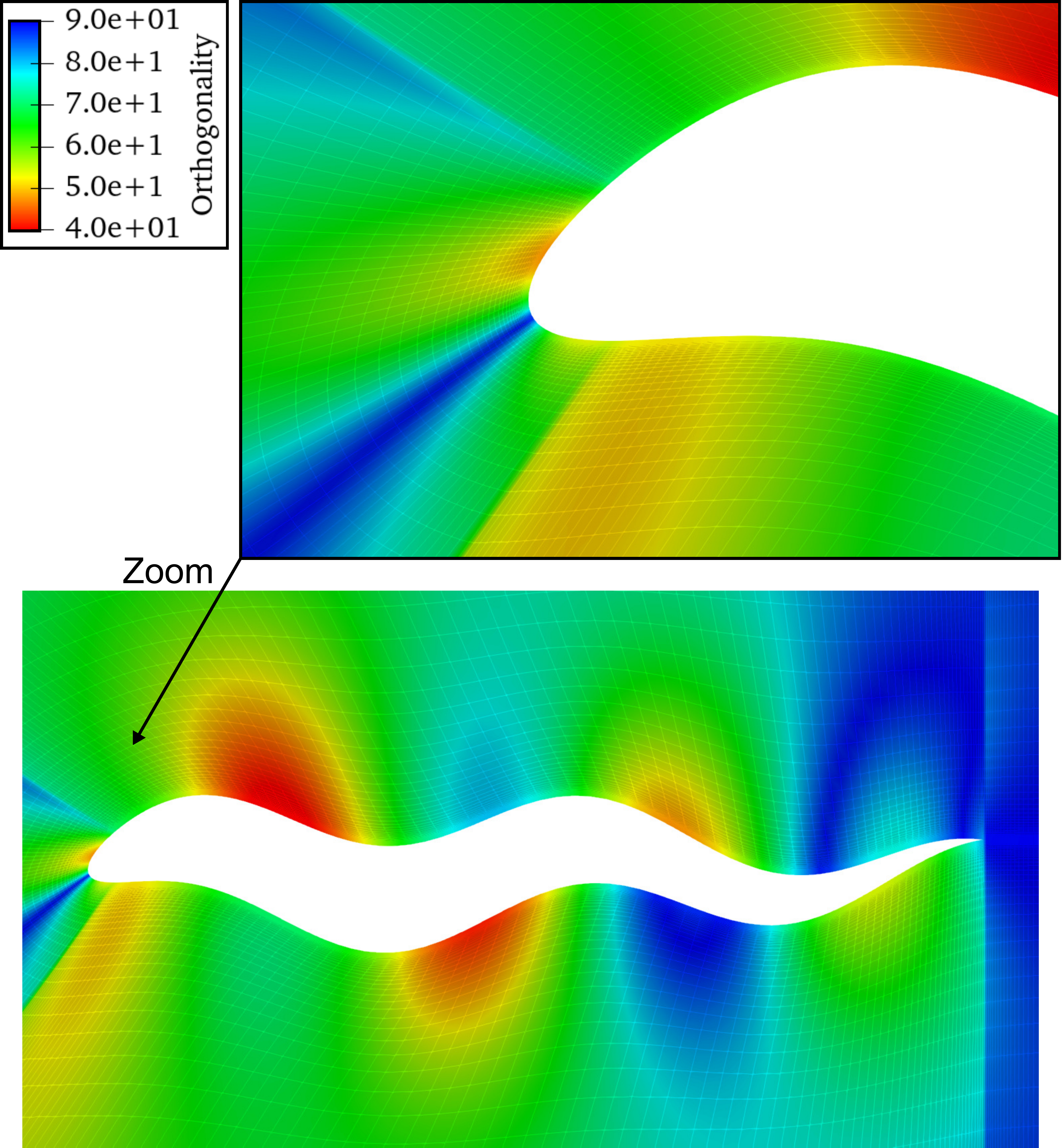}
}
\caption[Airfoil mesh quality]{Influence of the kind of deformation on the mesh quality. Comparing localised deformation to global deformation.}
\label{fig.7:airfoil_meshQuality}
\end{figure}

The influence of the kind of deformation on the mesh quality is shown in Fig.~\ref{fig.7:airfoil_meshQuality}. The mesh quality is evaluated by the orthogonality angle. The results show that a relatively high orthogonality angle is preserved throughout the simulations when compared to the clean mesh as shown in Fig.~\ref{fig.7:subfig-1:airfoil_meshQuality}. The localised ice deformation shown in Fig.~\ref{fig.7:subfig-2:airfoil_meshQuality} causes an isolated reduction in the mesh quality at the leading edge. While the global sinusoidal deformation shown in Fig.~\ref{fig.7:subfig-3:airfoil_meshQuality} causes a more universal reduction in the mesh quality with the greatest degradation located at the peaks of the waves.

\subsection{3D -- Swept Wing}
The aim of the second test case is to assess the efficiency and robustness of radial basis function mesh deformation data-reduction schemes on a three-dimensional icing problem. Experimental icing tests on a swept-wing with a NACA0012 profile performed in the NASA Glenn Icing Research Tunnel (IRT) \cite{tsao2012evaluation} are used. The results presented are part of an effort to help develop appropriate scaling methods for swept-wing ice accretion due to there being distinct morphological differences on straight and swept wings. One subset of these conditions is chosen and is outlined in Table~\ref{tab:3}. The test was carried out on the $15^{\textrm{th}}$ April, 2010 and corresponded to run number 07. The sweep angle of the wing was set at $45^{\circ}$ and the wing has a chord length of $0.914$ m. All tests were conducted at $0^{\circ}$ angle-of-attack. The freestream velocity is $51.44$ m/s. The outside air temperature is $260.5$ K The mean volume diameter of the supercooled water droplets is $44$ $\mu \textrm{m}$ which contain a liquid water content of $0.57$ $\textrm{g}/\textrm{m}^{3}$.   

\begin{figure}[htb!]
\centering
\subfloat[Unstructured volume mesh in the $x-z$ cut-plane at $y = 0$.\label{fig.8:subfig-1:wing_mesh}]
{
\includegraphics[width=0.48\linewidth]{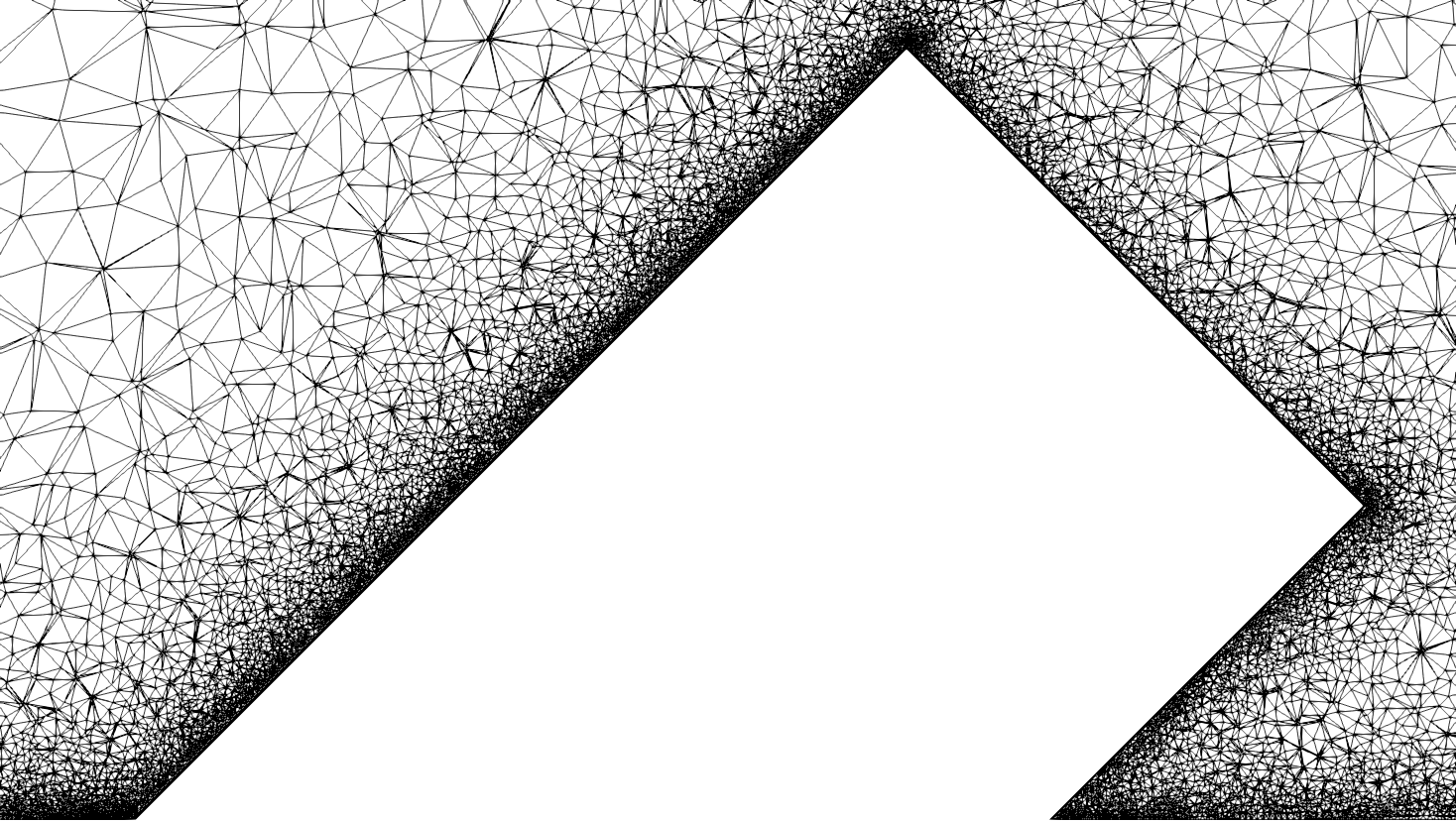}
} \hfill
\subfloat[Close-up view of the leading edge volume mesh in the $x-y$ cut-plane at $ z = 0.5 b$.\label{fig.8:subfig-2:wing_mesh}]
{
\includegraphics[width=0.48\linewidth]{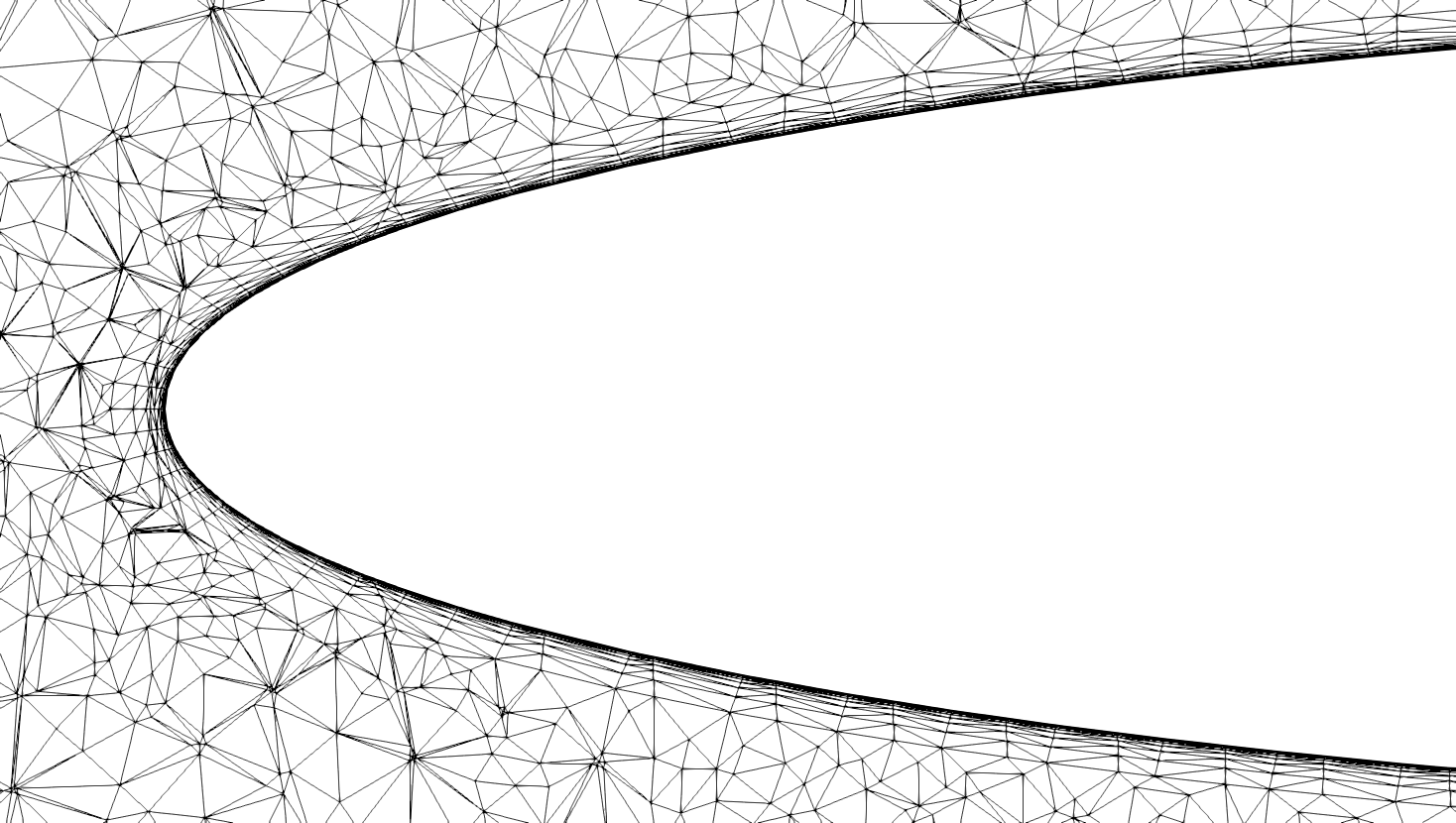}
}
\caption[Unstructured swept wing mesh]{Unstructured swept wing mesh. Constructed using a NACA0012 airfoil and based on a $45^{\circ}$ sweep angle.}
\label{fig.8:wing_mesh}
\end{figure}

The spatial discretization is achieved using an unstructured mesh which is shown in Fig.~\ref{fig.8:wing_mesh}. The total number of elements in the domain is $2,030,599$ and the total number of vertices is $591,225$. The swept wing has $21,650$ vertices distributed around its surface with there being an even distribution around the wing. Quadrilateral surface elements are used entirely along the leading edge. Elsewhere, to aid the generation of a high-quality surface mesh, quadrilateral elements dominate with the occasional triangular element. The far-field is placed 20 chord lengths from the wing. The grid resolution is sufficient at the wall to ensure $\textrm{y}+ < 1$.

\begin{figure}[htb!]
\centering
\includegraphics[width=0.98\linewidth]{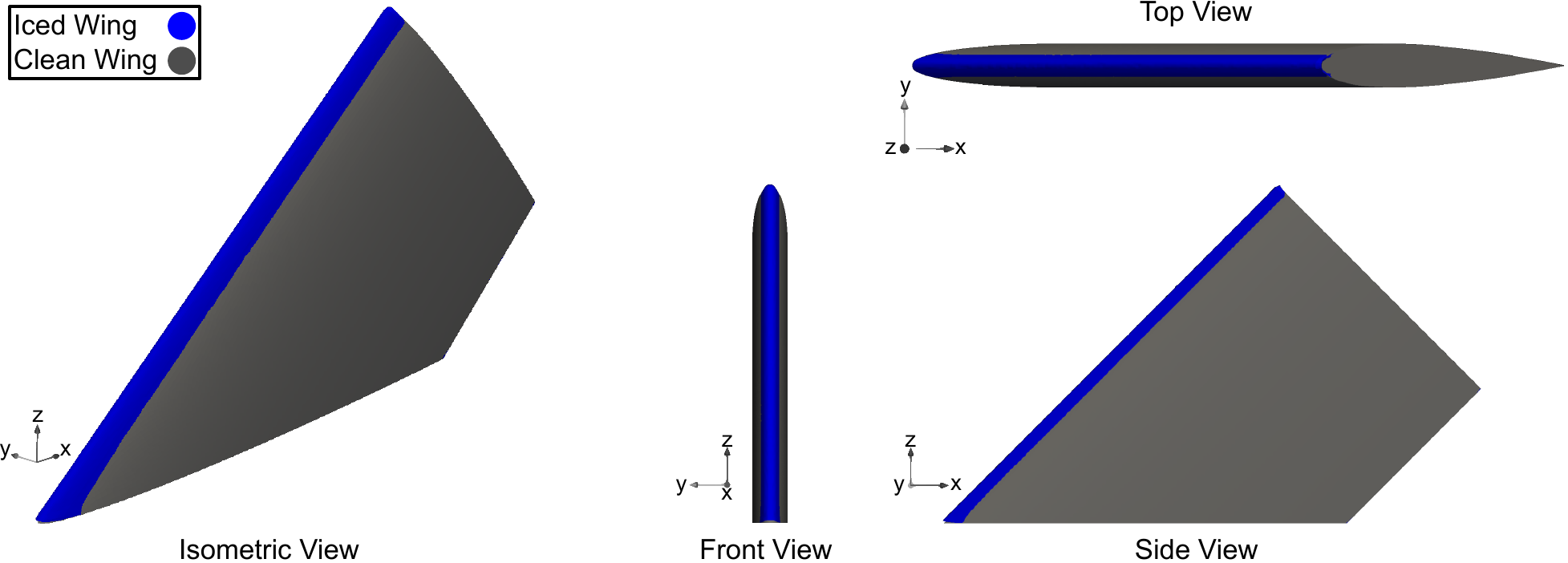}
\caption[Predicted ice shape on the swept wing]{Predicted 3D ice shape on the swept wing under conditions identified in Table~\ref{tab:3}.}
\label{fig.9:wing_ice}
\end{figure}

The SU2 solver configuration settings used to determine the aerodynamic flow field of the two-dimensional case are also used for the three-dimensional swept-wing case. The same applies to the ice accretion simulation settings of PoliMIce. For this case the ice prediction is computed using a single-step process. The final predicted ice shape after 19.8 minutes of ice accretion is shown in Fig.~\ref{fig.9:wing_ice}. The ice shape exhibits a spearhead like shape paradigmatic of the rime regime along the leading edge of the wing.

The radial basis function mesh deformation uses the compactly supported Wendland C2 function. Points within the support radius of $R = 3c$ are deformed. Three levels of multi-level greedy surface point reduction are used. Each level is updated when there is a reduction of $\varepsilon = 10^{-1}$ in the normalized error. A volume reduction factor of $k = 5$ is chosen. 

\begin{figure}[htb!]
\centering
\includegraphics[width=0.98\linewidth]{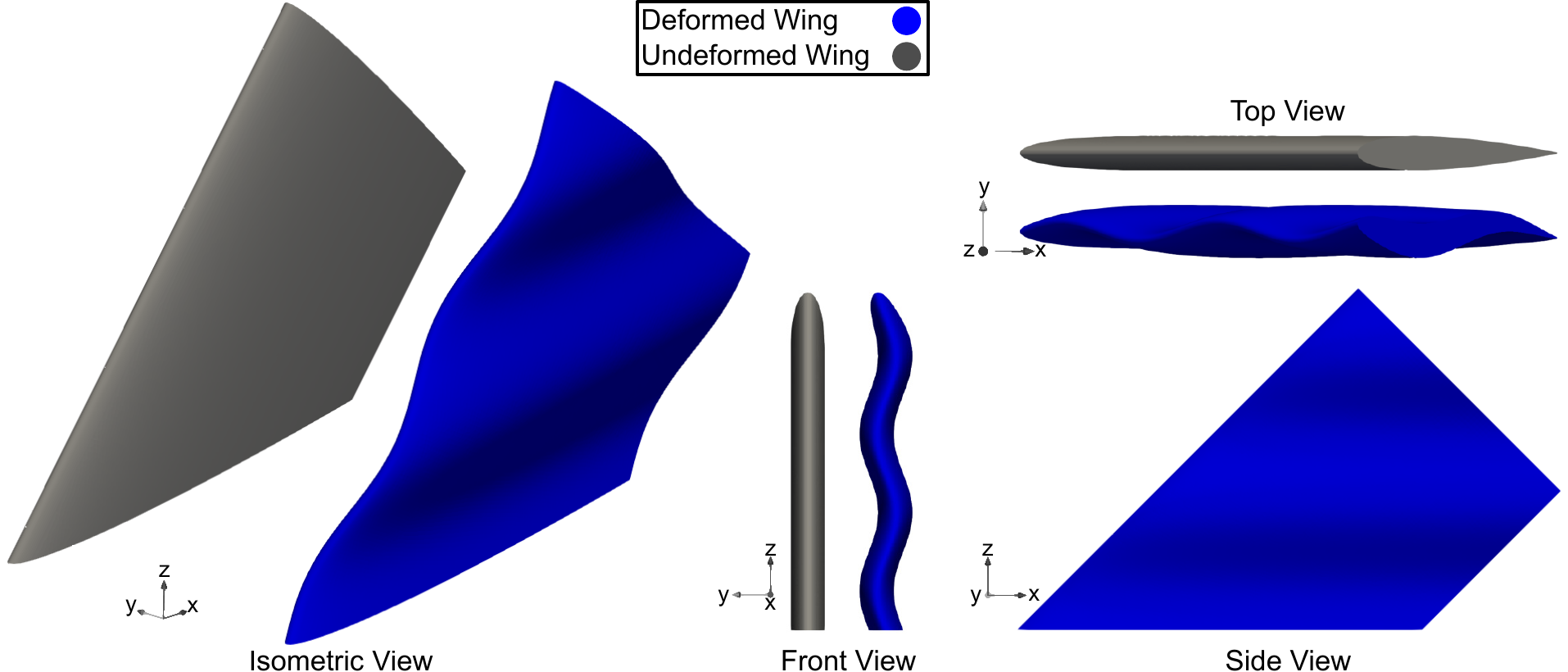}
\caption[Sinusoidal motion of the wing]{Sinusoidal motion of the wing described by Eq.~\ref{eq:27}.}
\label{fig.10:wing_sinusiodal}
\end{figure}

The two-dimensional benchmark test case is extended to evaluate the performance of the three-dimensional swept-wing test case. Sinusoidal motion is thus applied to the swept-wing and can be described by:

\begin{equation}
\label{eq:27}
	\Delta y = 0.03 \sin \left ( 4 \pi z \right ) \, ,
\end{equation}
where $\Delta y$ describes the displacement of the swept-wing as a function of its position in the $z$-direction along its span, $b$. The deformation applied to the swept-wing is shown in Fig.~\ref{fig.10:wing_sinusiodal}. The mesh and radial basis function parameters are again maintained homogeneous to ensure maximum similarity between the localised iced deformation and the global sinusoidal deformation. 

\begin{figure}[htb!]
\centering
\subfloat[Convergence history in terms of selected points.\label{fig.11:subfig-1:wing_greedyConvergence}]
{
\includegraphics[width=0.6\linewidth]{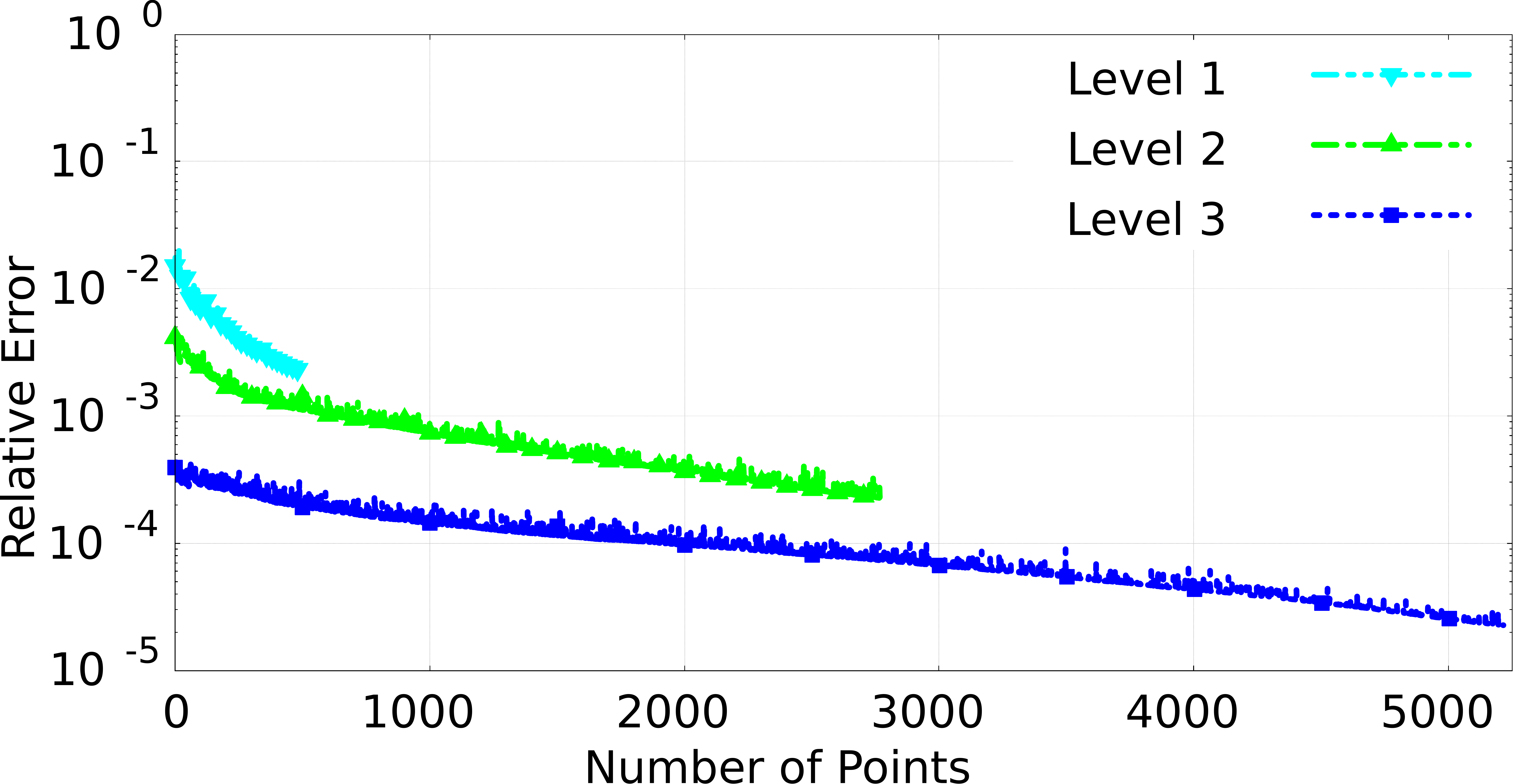}
} \\
\subfloat[Convergence history in terms of CPU time.\label{fig.11:subfig-2:wing_greedyConvergence}]
{
\includegraphics[width=0.6\linewidth]{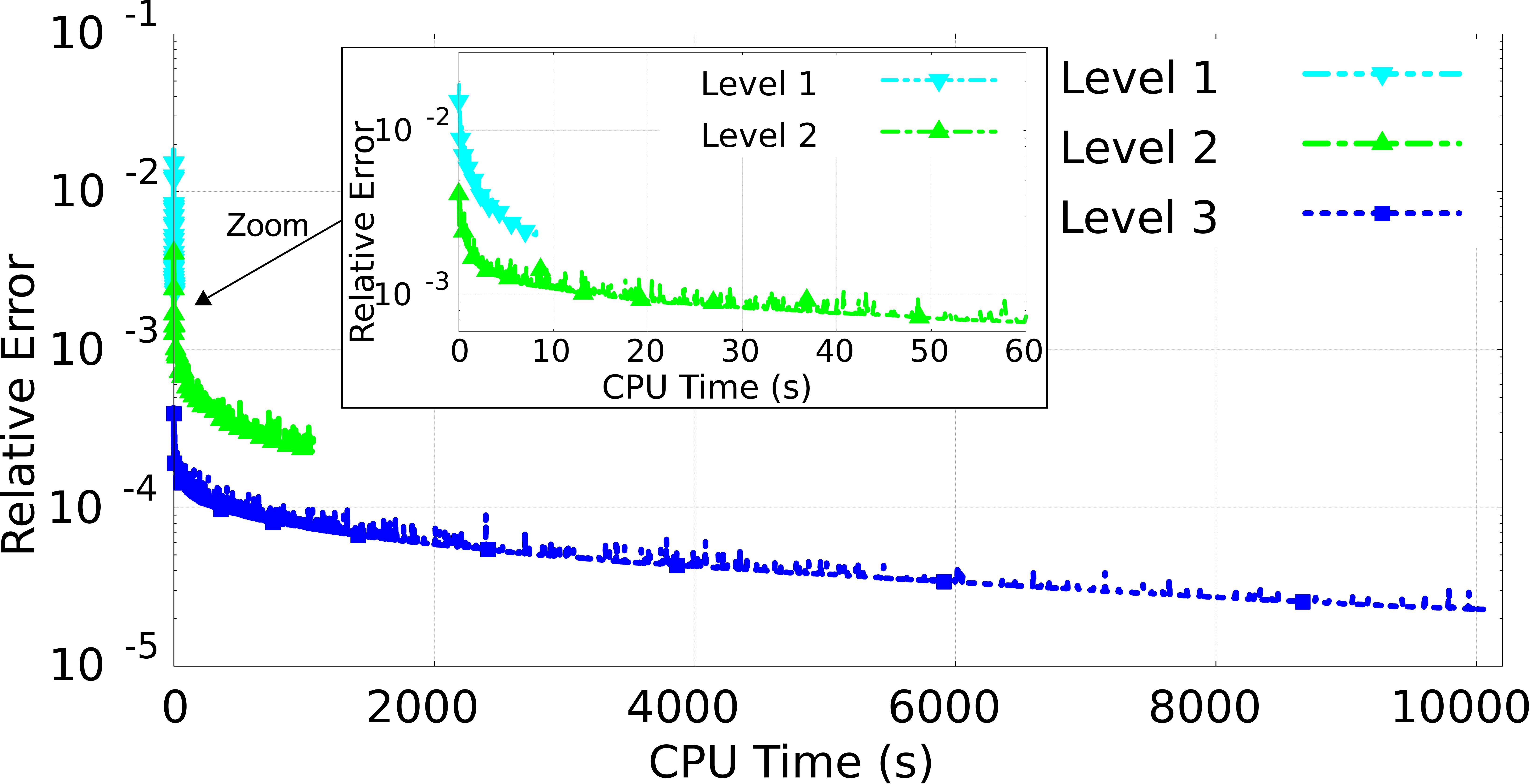}
}
\caption[Comparison of error reduction rates on the iced wing]{Comparison of error reduction rates in terms of selected points and CPU time for the NACA0012 swept wing post deformation due to icing.}
\label{fig.11:wing_greedyConvergence}
\end{figure}

\begin{figure}[htb!]
\centering
\subfloat[Convergence history in terms of selected points.\label{fig.12:subfig-1:defWing_greedyConvergence}]
{
\includegraphics[width=0.6\linewidth]{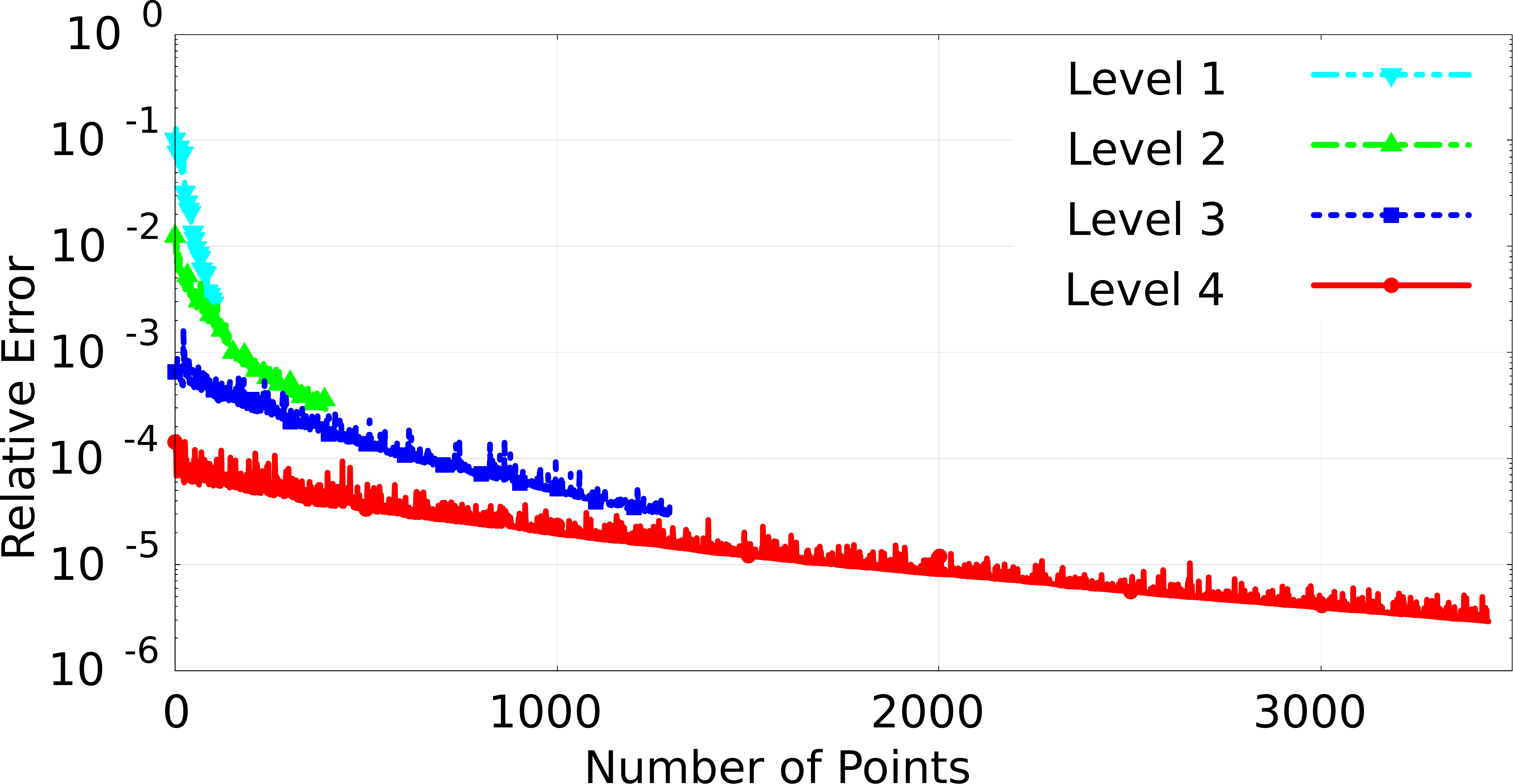}
} \\
\subfloat[Convergence history in terms of CPU time.\label{fig.12:subfig-2:defWing_greedyConvergence}]
{
\includegraphics[width=0.6\linewidth]{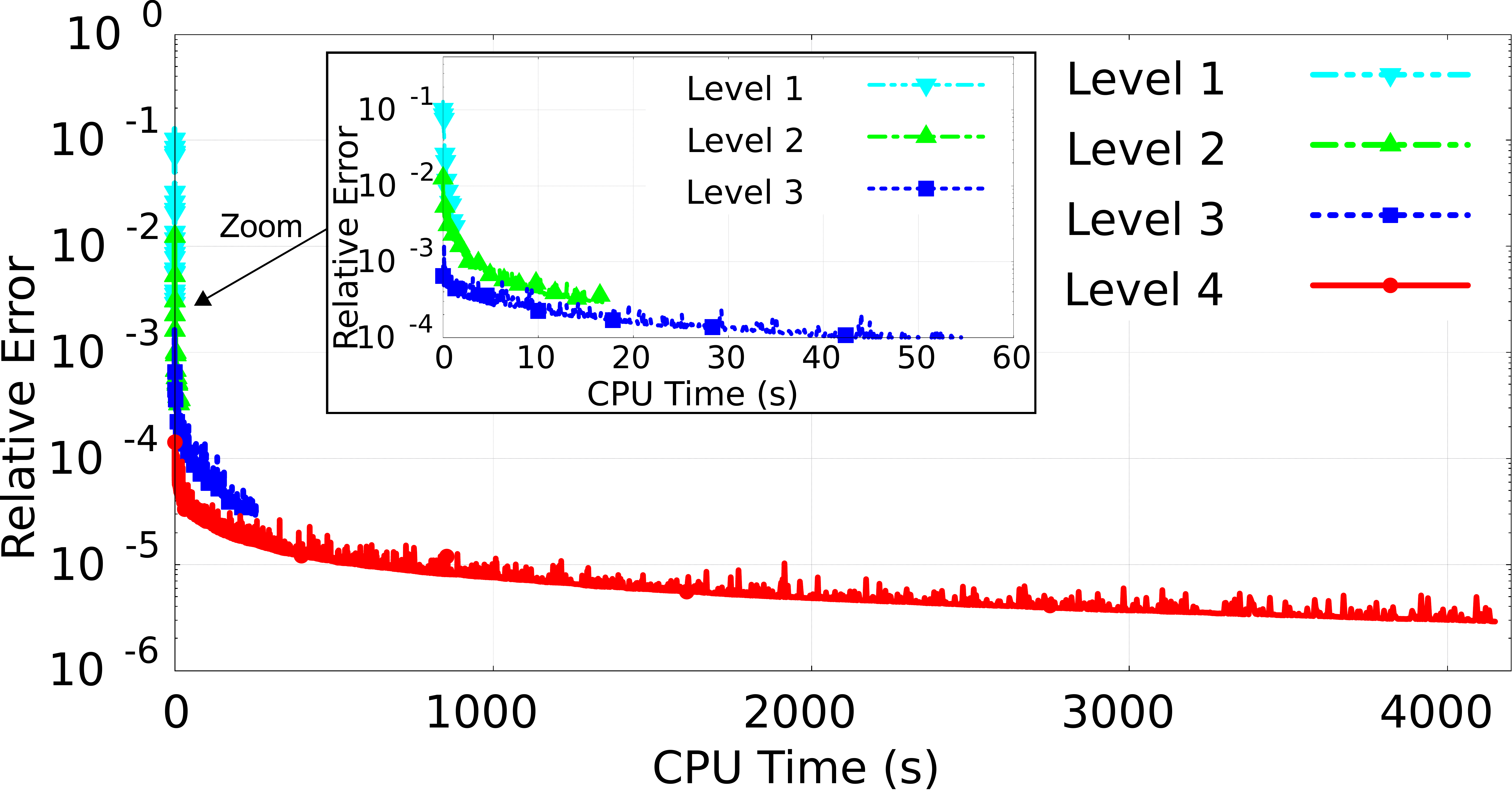}
}
\caption[Comparison of error reduction rates on the sinusoidal wing]{Comparison of error reduction rates in terms of selected points and CPU time for the NACA0012 swept wing post deformation due to sinusoidal motion.}
\label{fig.12:defWing_greedyConvergence}
\end{figure}

The convergence history of the normalized displacement error during the iced and sinusoidal deformation of the swept-wing are respectively illustrated in Figs.~\ref{fig.11:wing_greedyConvergence} \&~\ref{fig.12:defWing_greedyConvergence}. Given the increase in size of the data-set, the efficiency of the multi-level greedy point selection performs well in terms of maintaining a low number of control points and in terms of CPU time. When comparing these two test cases it is clear that the global sinusoidal deformation performs better. While the localised iced deformation shown in Fig.~\ref{fig.11:wing_greedyConvergence} can obtain three-levels of greedy surface point selection, the global sinusoidal deformation shown in Fig.~\ref{fig.12:defWing_greedyConvergence} can obtain four-levels of greedy surface point selection. In essence, while global sinusoidal deformations can be characterised by a small number of control points, localised iced deformations simply cannot. Despite this, satisfactory reduction in the normalised surface error is achieved by both test cases. Most significantly of all, the data-reduction techniques do indeed help to reduce the high computational cost associated to radial basis function mesh deformation on large data-sets. Notably within the first 60 seconds of the iced and sinusoidal deformations, normalized surface errors of $10^{-3}$ and $10^{-4}$ are respectively achieved as revealed in Figs.~\ref{fig.11:subfig-2:wing_greedyConvergence} \&\ref{fig.12:subfig-2:defWing_greedyConvergence}. 

\begin{figure}[htb!]
\centering
\subfloat[Level 1: 492 Control Points.\label{fig.13:subfig-1:wing_greedySelection}]
{
\includegraphics[width=0.48\linewidth]{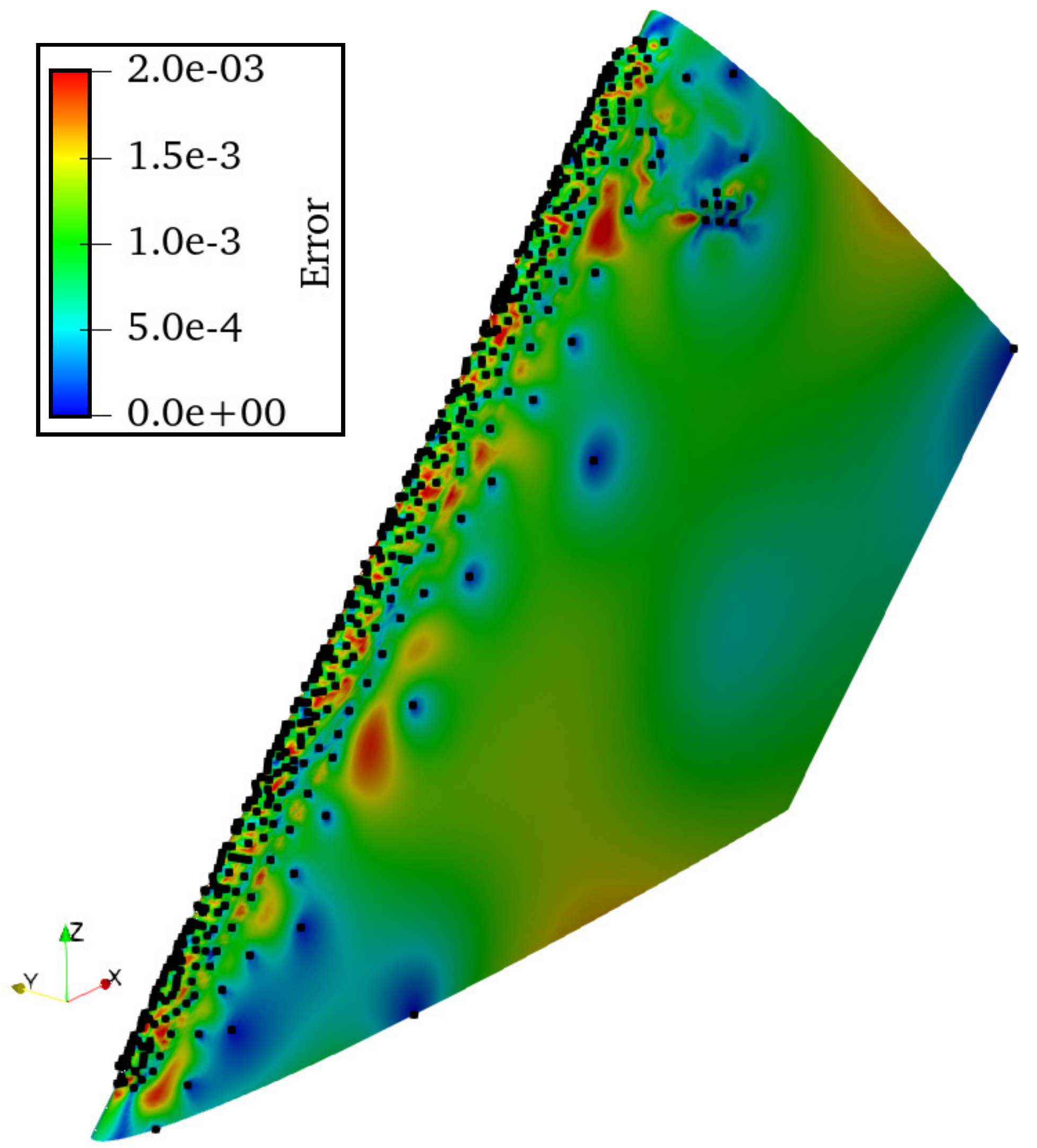}
} \hfill
\subfloat[Level 2: 2767 Control Points.\label{fig.13:subfig-2:wing_greedySelection}]
{
\includegraphics[width=0.48\linewidth]{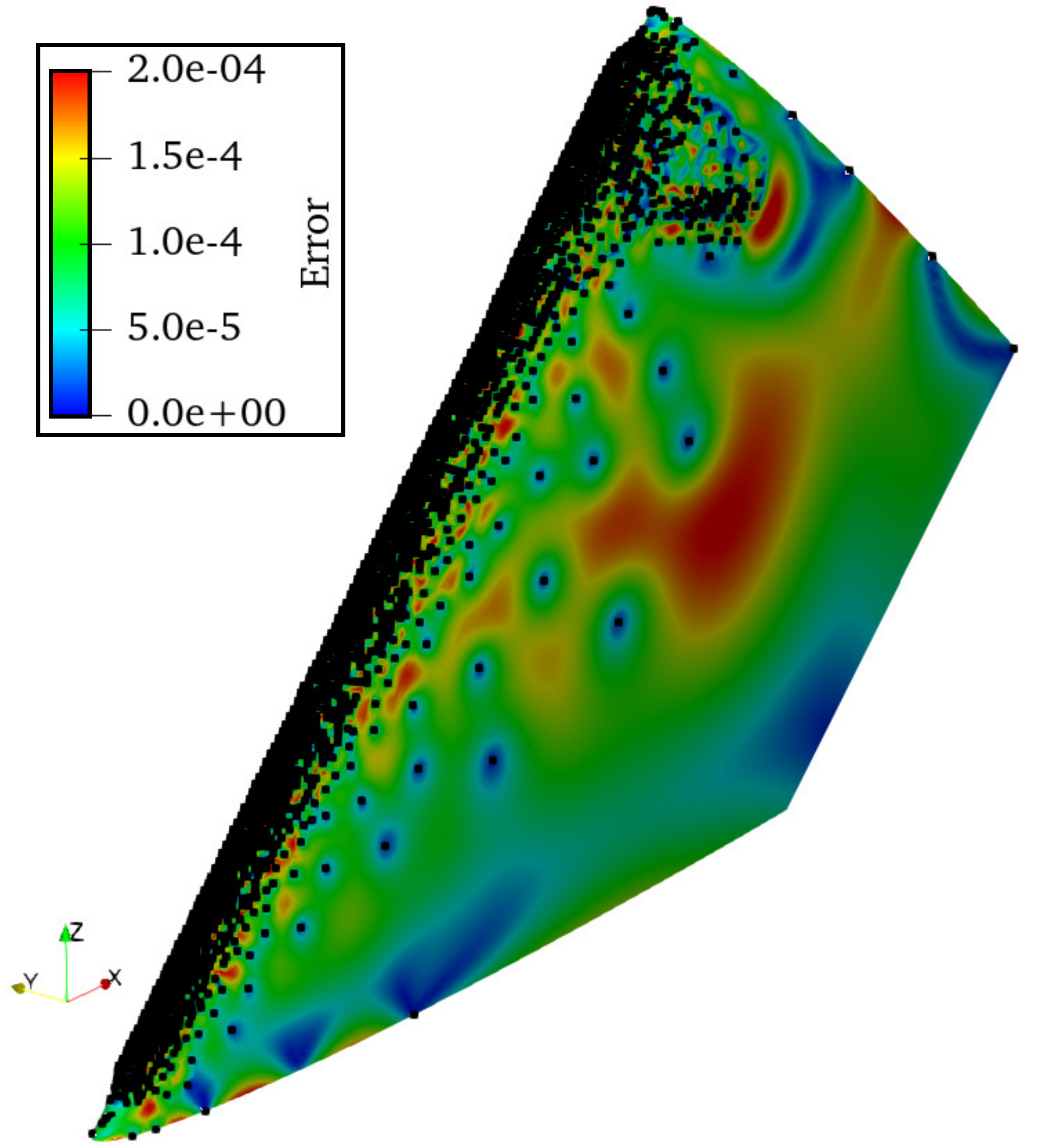}
} \\
\subfloat[Level 3: 5214 Control Points.\label{fig.13:subfig-3:wing_greedySelection}]
{
\includegraphics[width=0.48\linewidth]{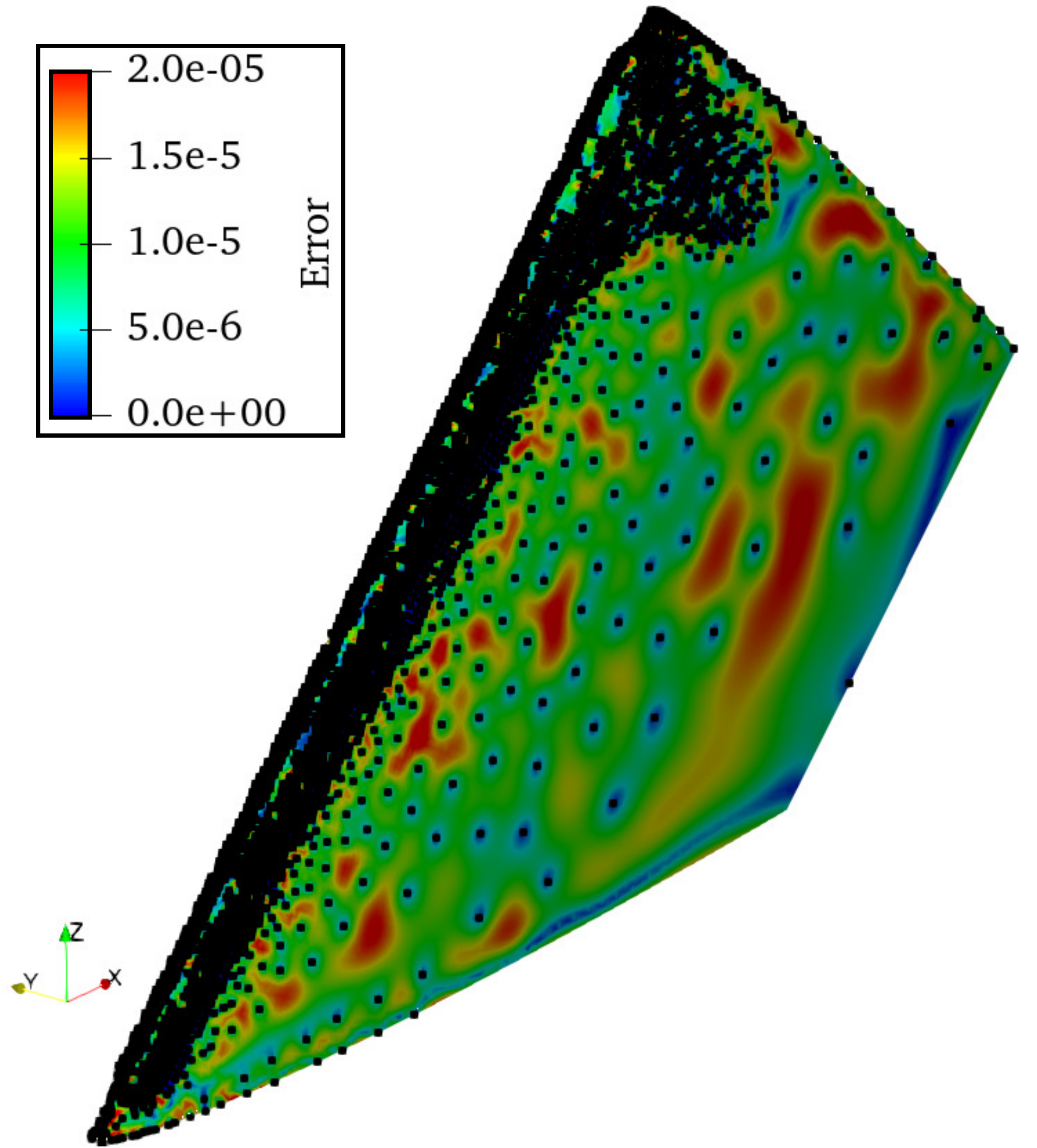}
} 
\caption[Control points selected on the iced wing]{Contour map of the normalized surface error during the multi-level greedy point selection on the iced wing. The black points indicate the control points.}
\label{fig.13:wing_greedySelection}
\end{figure}

\begin{figure}[htb!]
\centering
\subfloat[Level 1: 102 Control Points.\label{fig.14:subfig-1:defWing_greedySelection}]
{
\includegraphics[width=0.48\linewidth]{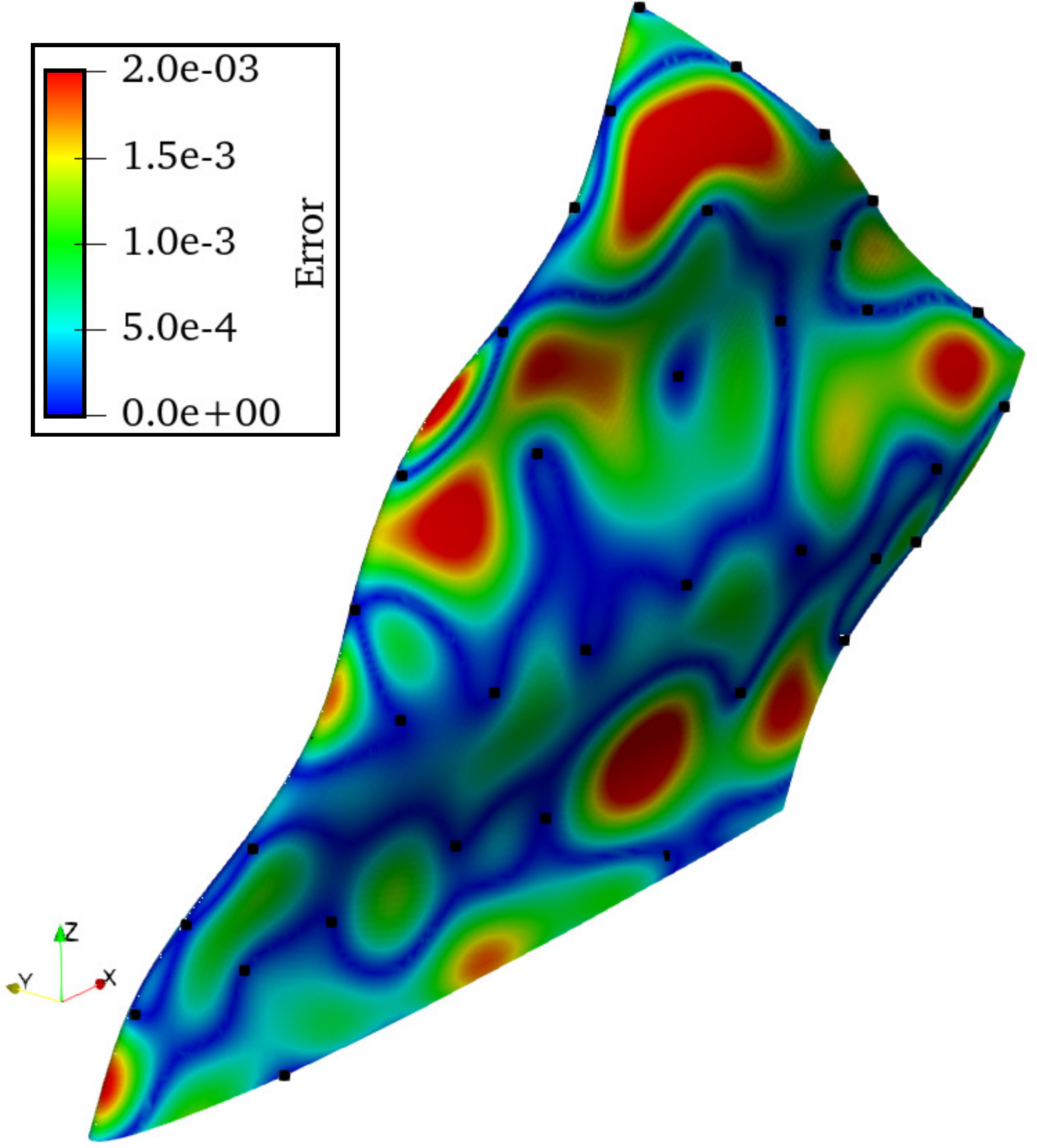}
} \hfill
\subfloat[Level 2: 395 Control Points.\label{fig.14:subfig-2:defWing_greedySelection}]
{
\includegraphics[width=0.48\linewidth]{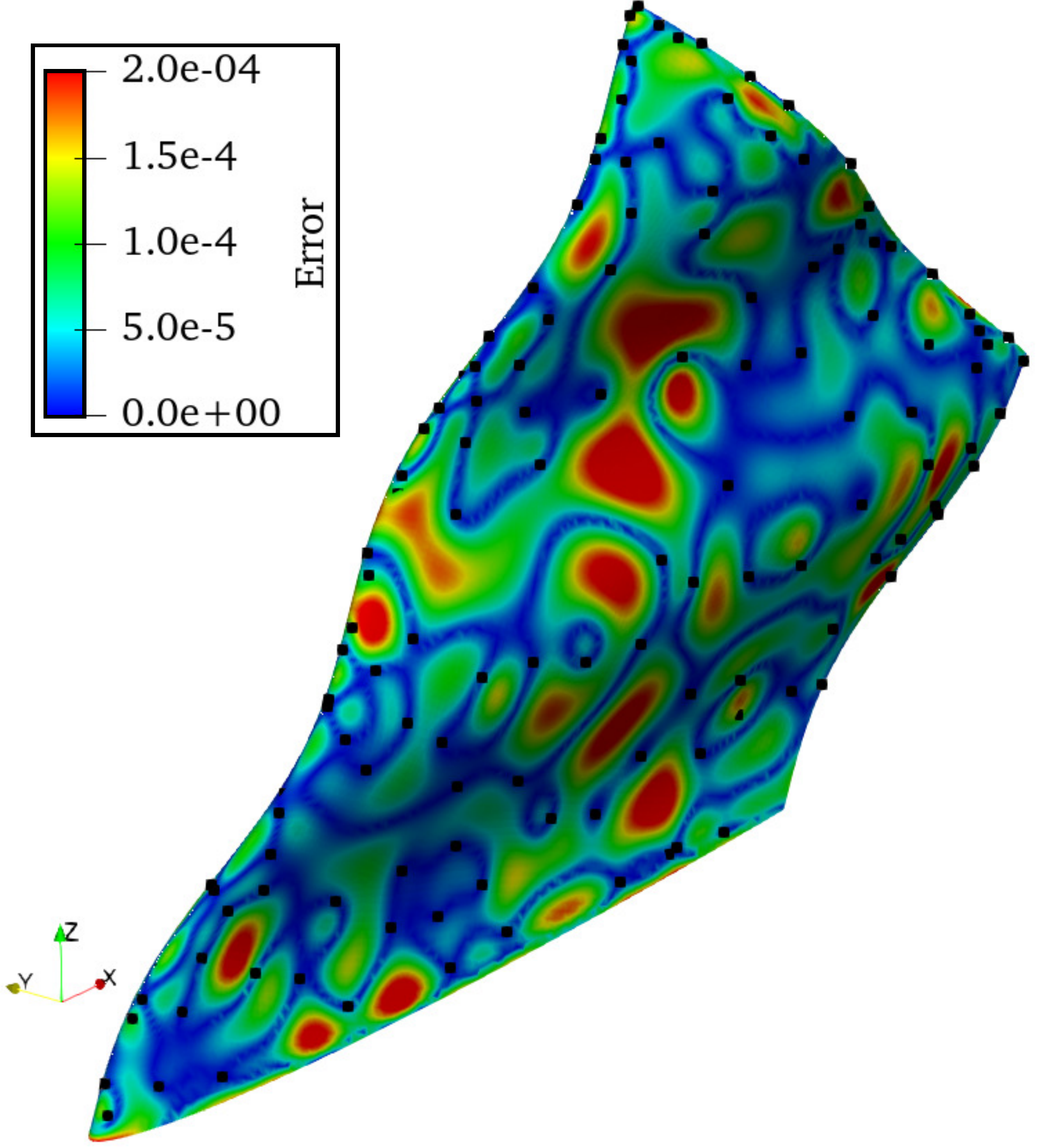}
} \\
\subfloat[Level 3: 1297 Control Points.\label{fig.14:subfig-3:defWing_greedySelection}]
{
\includegraphics[width=0.48\linewidth]{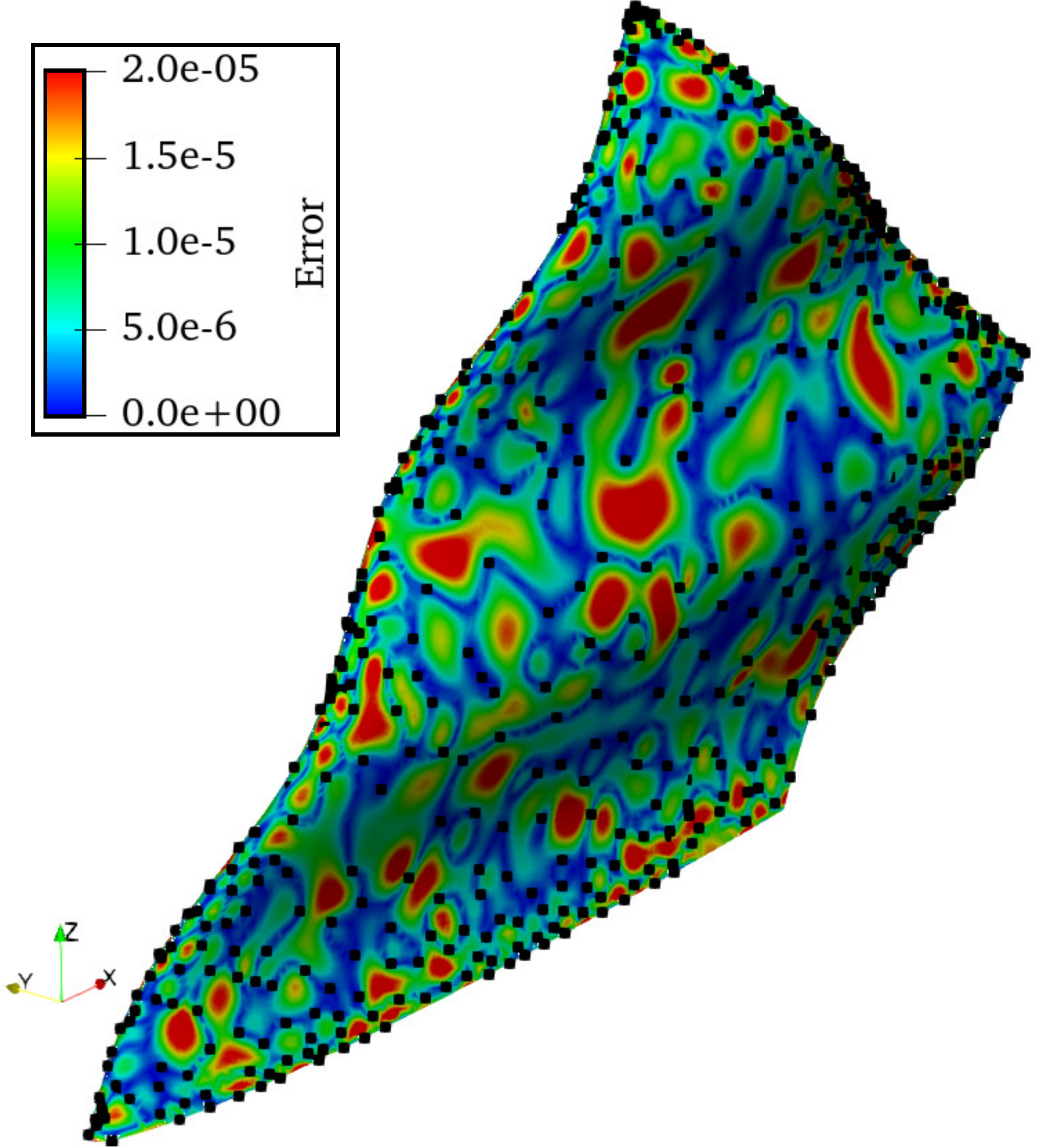}
} \hfill
\subfloat[Level 4: 3437 Control Points.\label{fig.14:subfig-4:defWing_greedySelection}]
{
\includegraphics[width=0.48\linewidth]{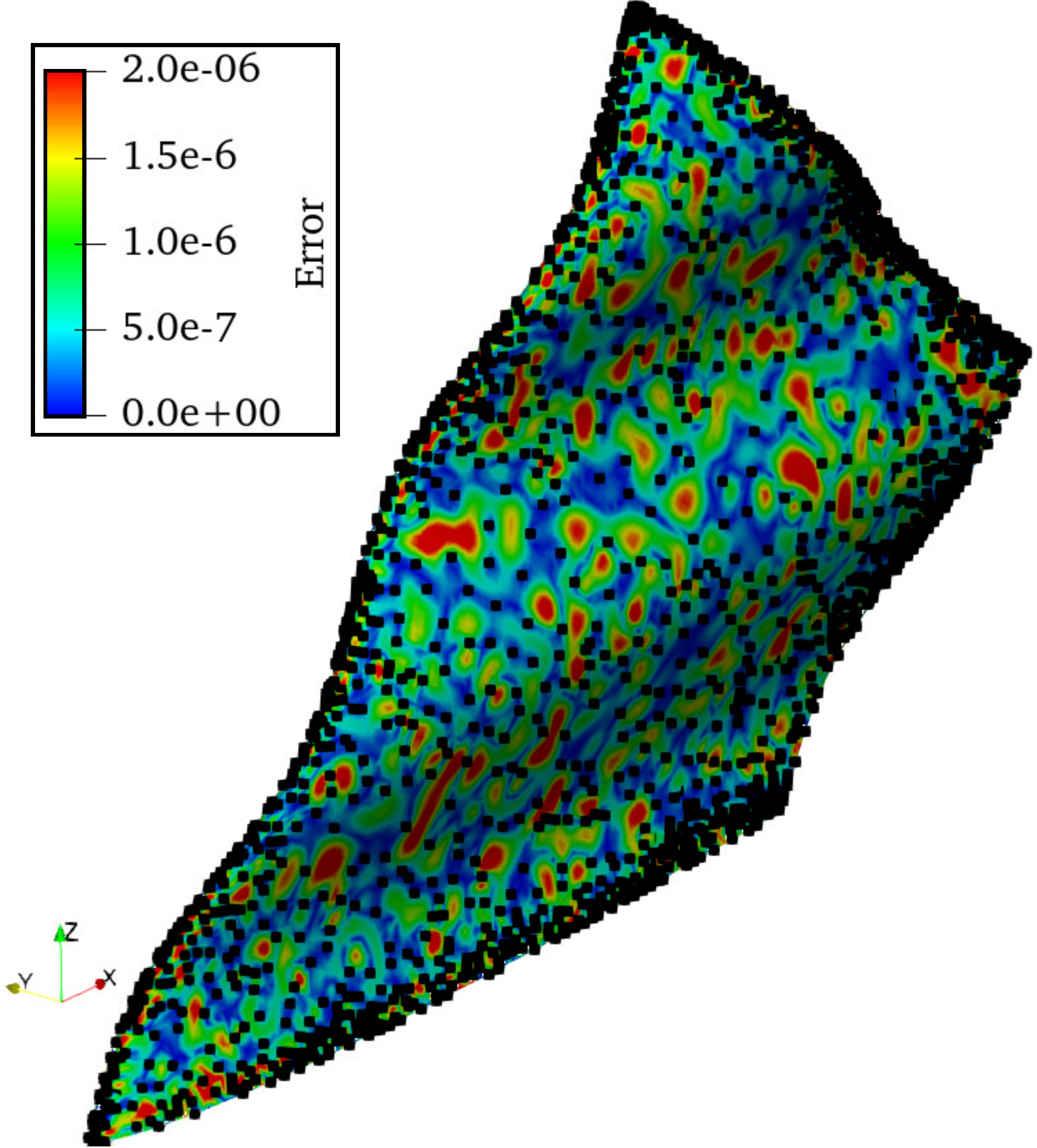}
} 
\caption[Control points selected on the sinusoidal wing]{Contour map of the normalized surface error during the multi-level greedy point selection on the sinusoidal. The black points indicate the control points.}
\label{fig.14:defWing_greedySelection}
\end{figure}

The control points selected throughout the multi-level greedy selection process during the iced and sinusoidal deformation are respectively shown in Figs.~\ref{fig.13:wing_greedySelection} \&\ref{fig.14:defWing_greedySelection}. The test cases exhibit significantly different distributions of control points. The local iced deformation illustrated in Fig.~\ref{fig.13:wing_greedySelection} depicts the intelligence of the selection process to use control points primarily congregated around the leading edge of the swept-wing where ice is accreted. Resultantly an anisotropic control point distribution is present. The global sinusoidal deformation visible in Fig.~\ref{fig.14:defWing_greedySelection} displays a more expansive distribution of selected control points around the swept-wing which reflects the deformation described by Eq.~\ref{eq:27}. 

\begin{figure}[htb!]
\centering
\subfloat[Orthogonality of the clean wing mesh prior to deformation.\label{fig.15:subfig-1:wing_meshQuality}]
{
\includegraphics[width=0.7\linewidth]{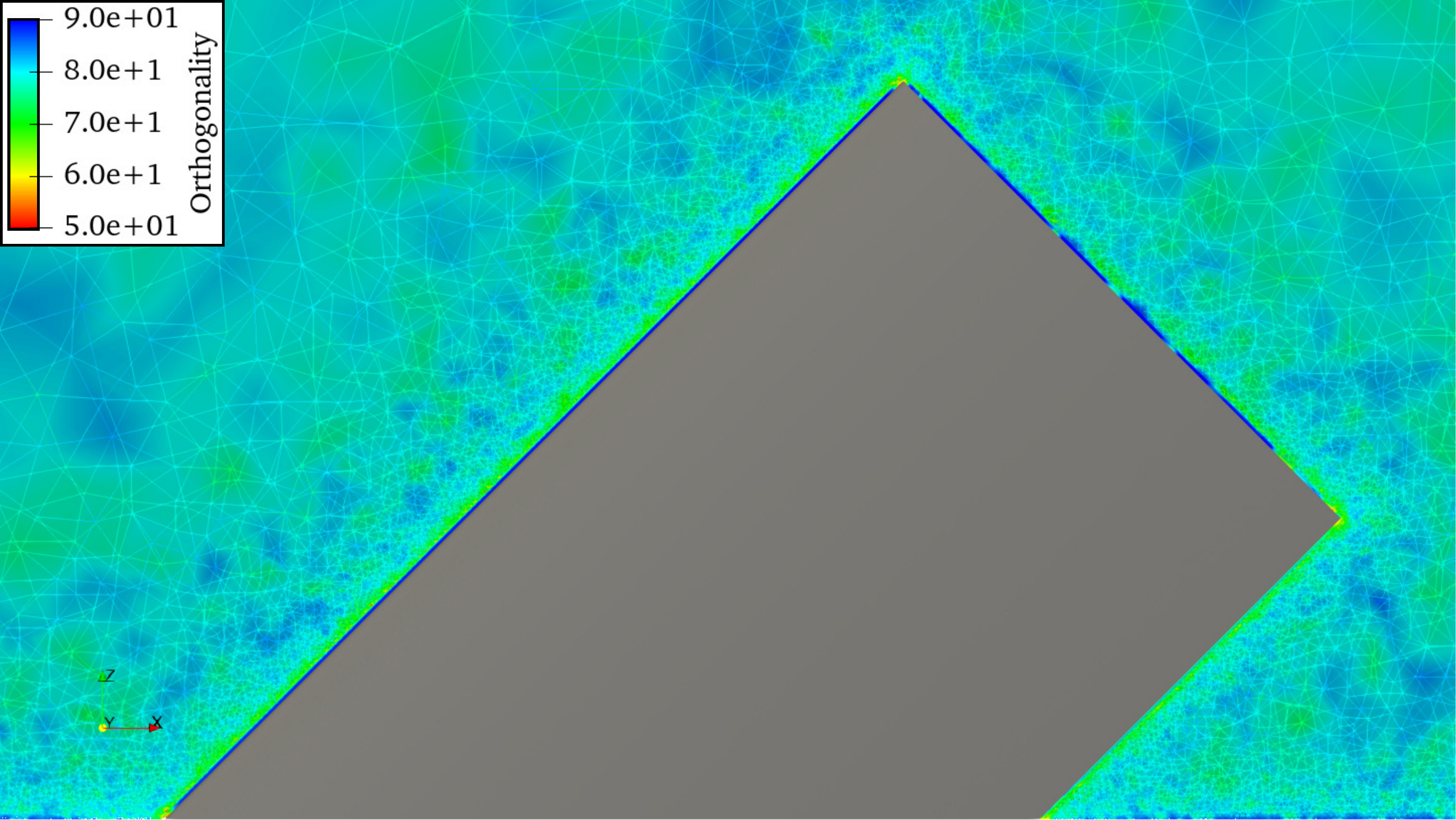}
} \\
\subfloat[Orthogonality of the iced wing mesh post deformation.\label{fig.15:subfig-2:wing_meshQuality}]
{
\includegraphics[width=0.7\linewidth]{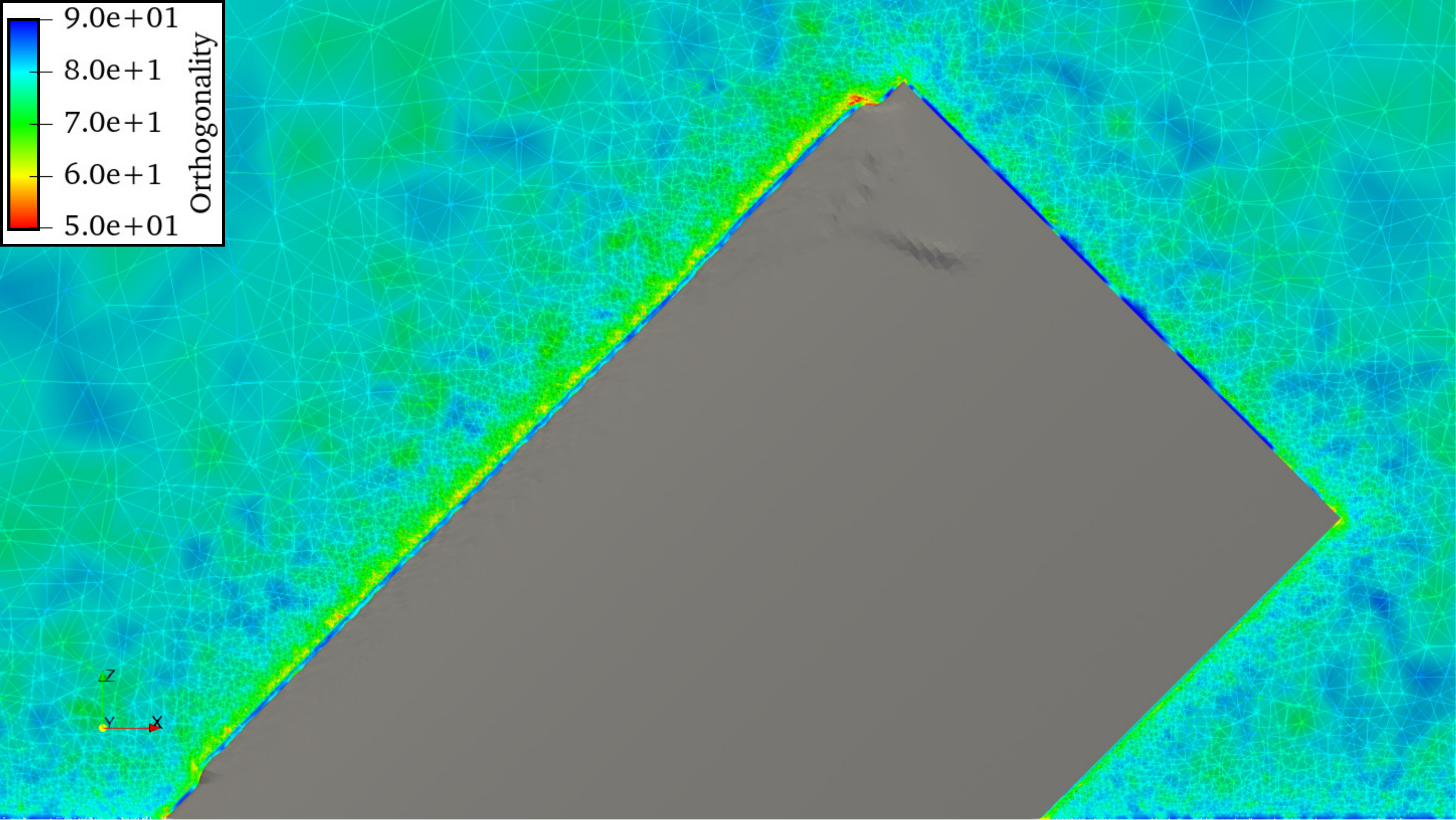}
} \\
\subfloat[Orthogonality of the wing mesh post sinusoidal deformation.\label{fig.15:subfig-3:wing_meshQuality}]
{
\includegraphics[width=0.7\linewidth]{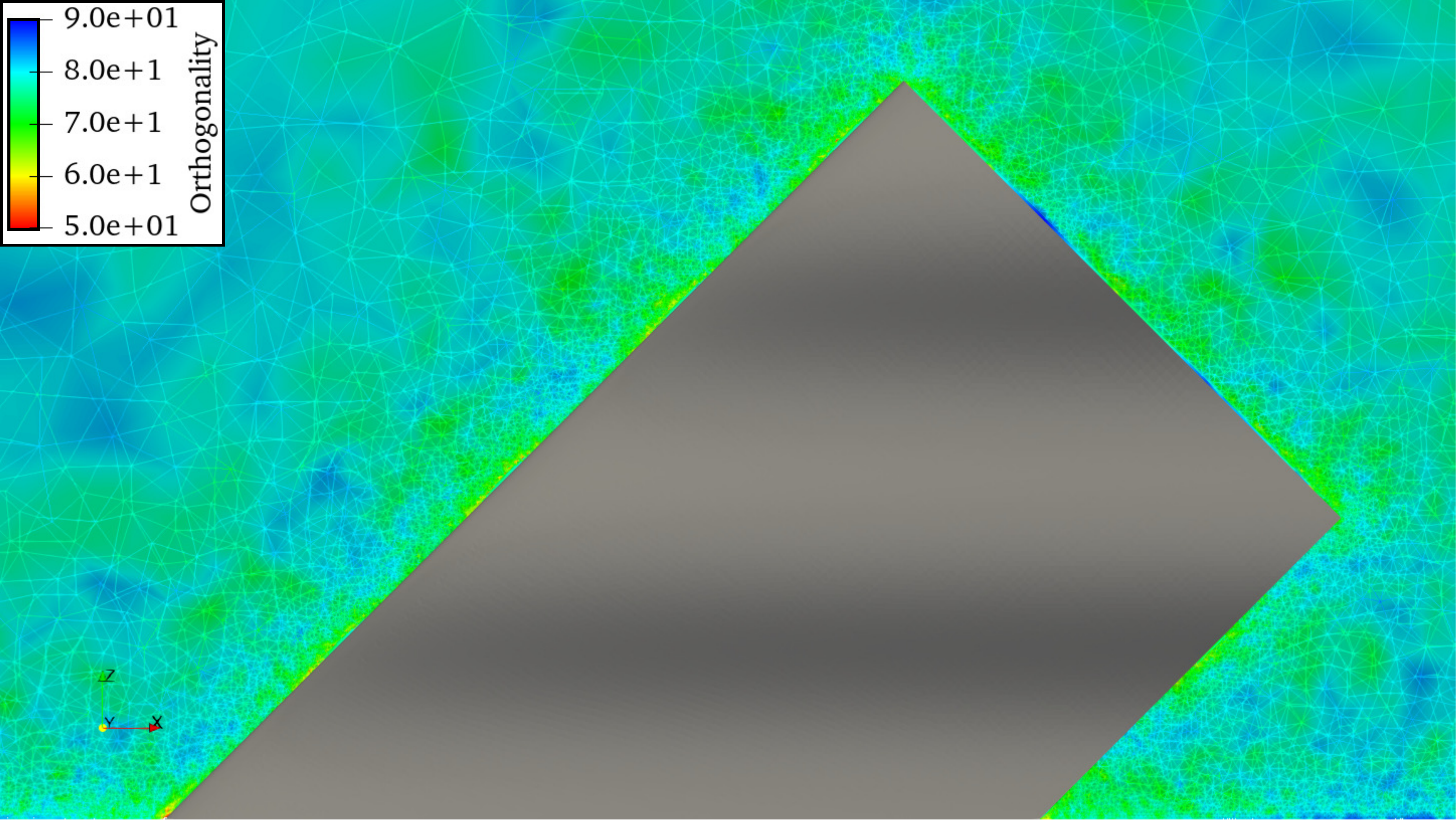}
}
\caption[Wing mesh quality in the $x-z$ plane]{Influence of the kind of mesh deformation on the quality of the swept wing mesh in the $x-z$ plane at $y = 0$.}
\label{fig.15:wing_meshQuality}
\end{figure}

\begin{figure}[htb!]
\centering
\subfloat[Orthogonality of the clean wing mesh prior to deformation.\label{fig.16:subfig-1:wing_meshQuality2}]
{
\includegraphics[width=0.4\linewidth]{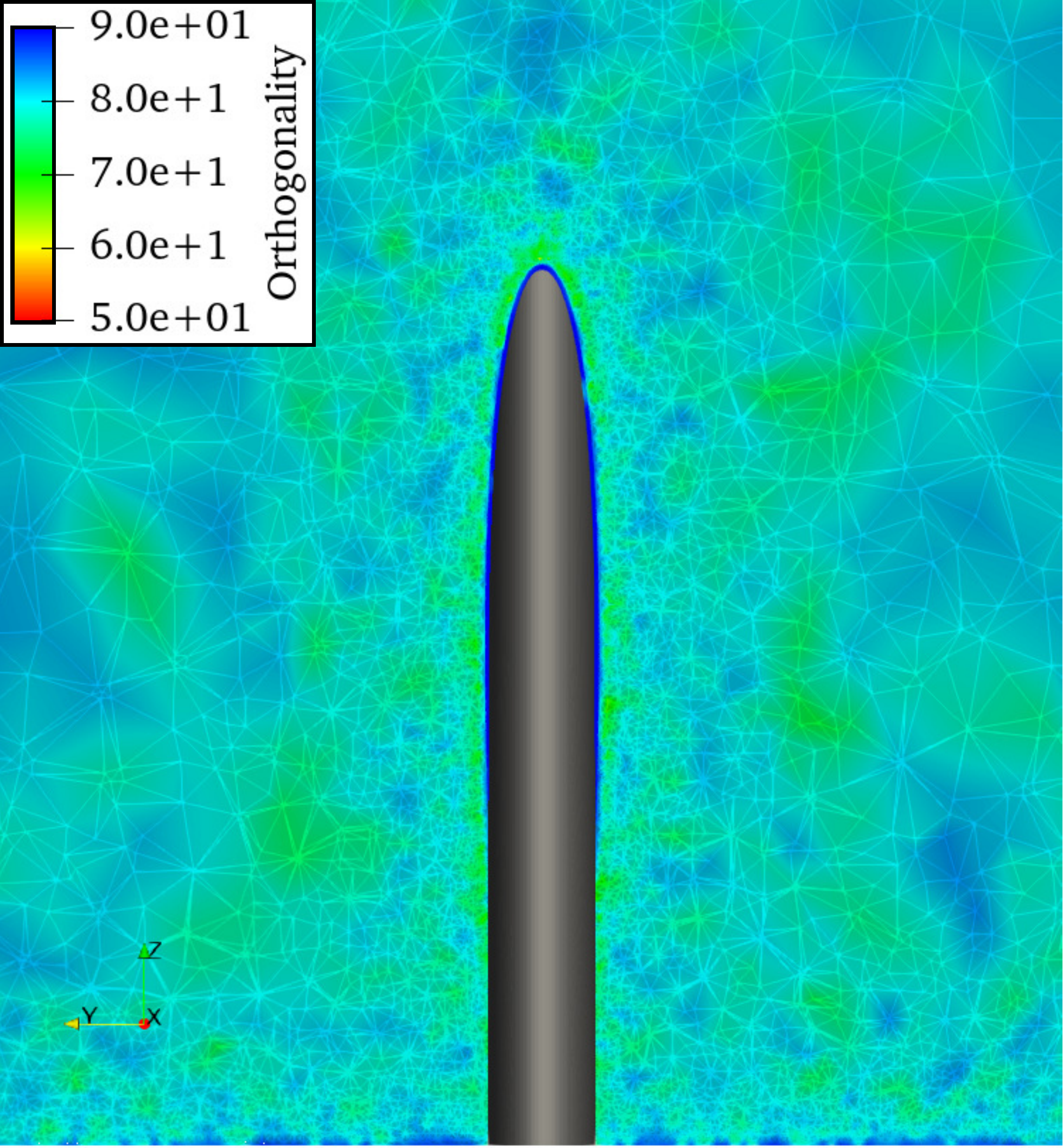}
} \hfill
\subfloat[Orthogonality of the iced wing mesh post deformation.\label{fig.16:subfig-2:wing_meshQuality2}]
{
\includegraphics[width=0.4\linewidth]{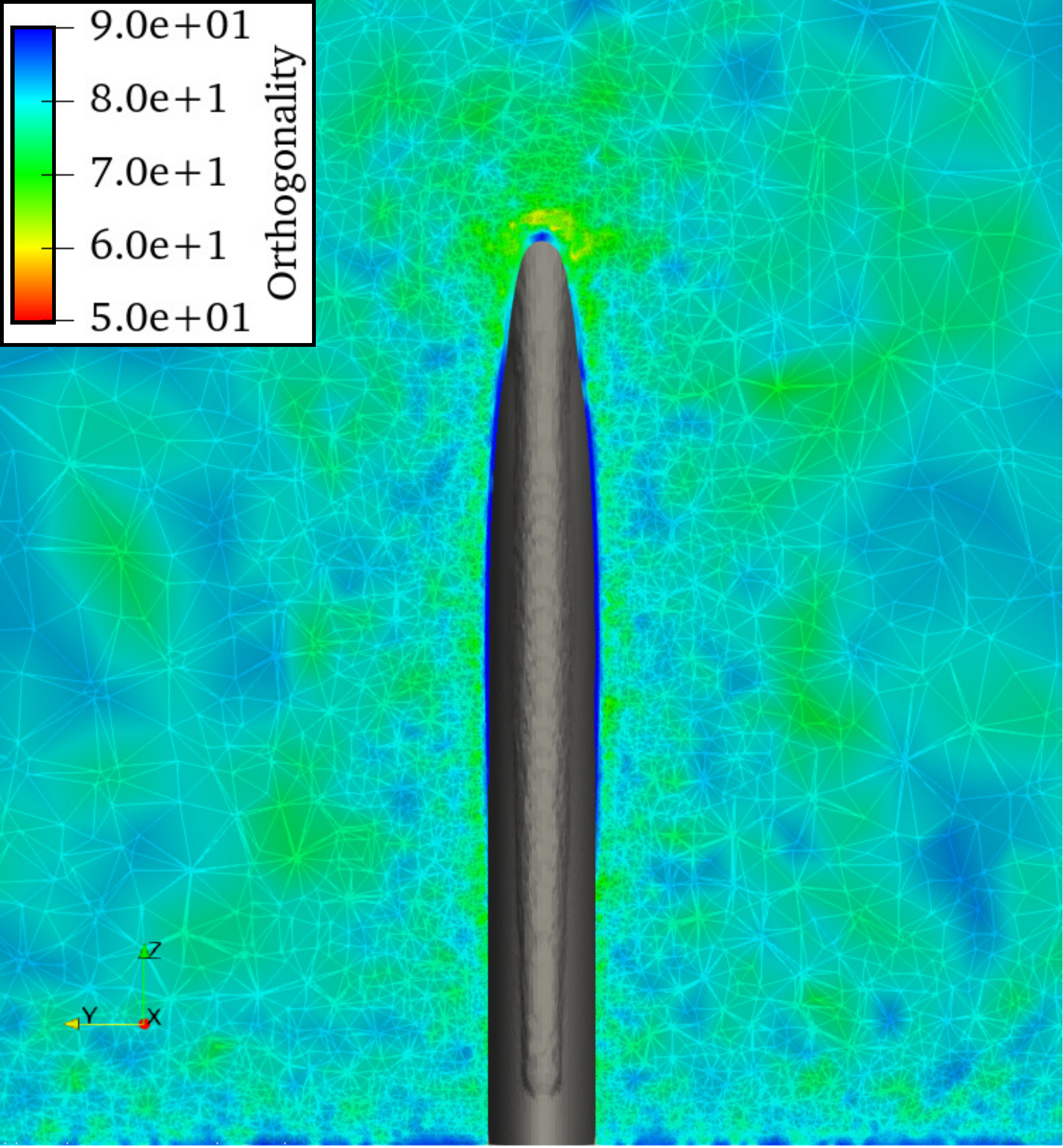}
} \\
\subfloat[Orthogonality of the wing mesh post sinusoidal deformation.\label{fig.16:subfig-3:wing_meshQuality2}]
{
\includegraphics[width=0.4\linewidth]{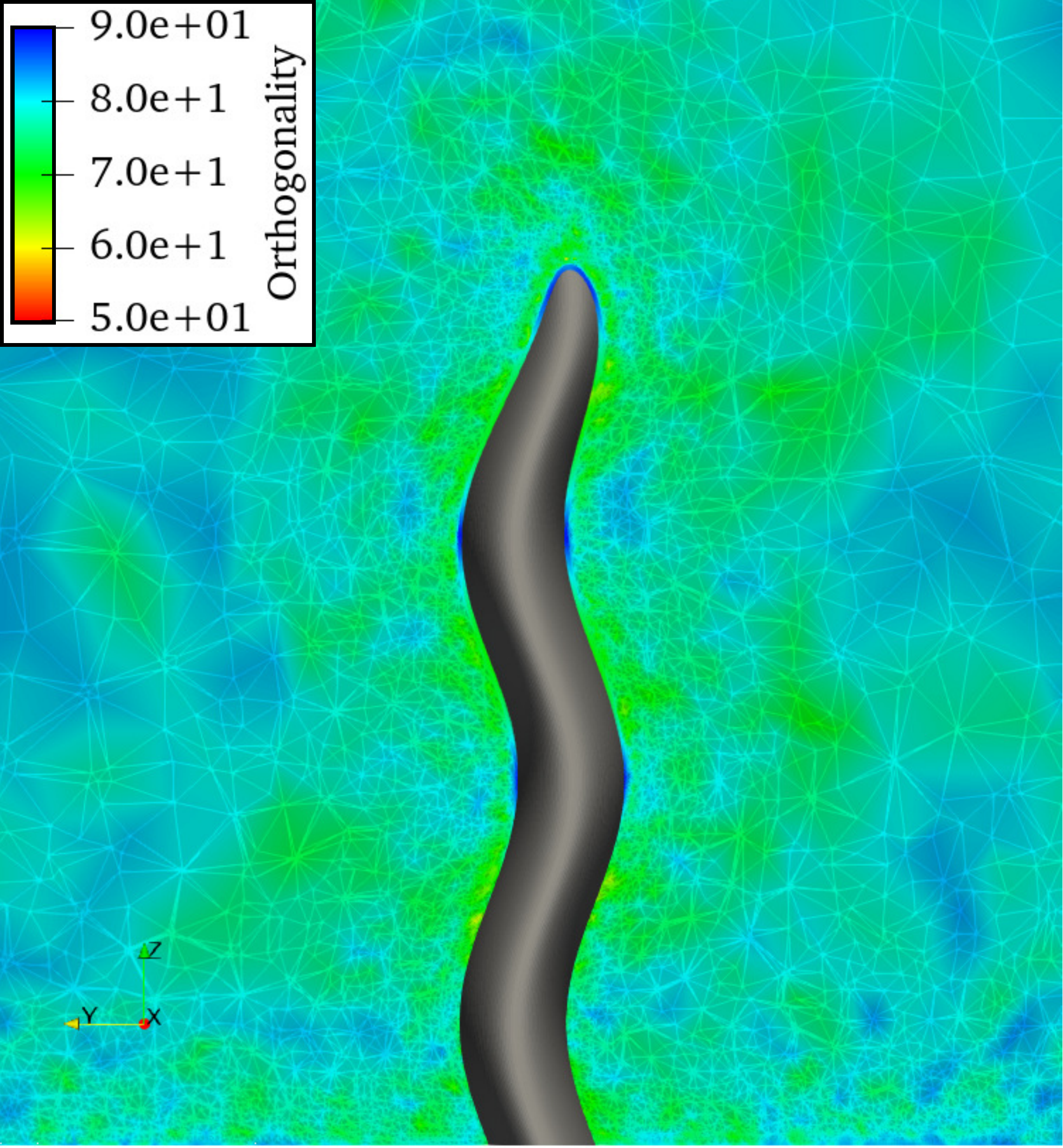}
}
\caption[Wing mesh quality in the $y-z$ plane]{Influence of the kind of mesh deformation on the quality of the swept wing mesh in the $y-z$ plane at $x = c$.}
\label{fig.16:wing_meshQuality2}
\end{figure}

The influence of the kind deformation on the mesh quality is shown in Figs.~\ref{fig.15:wing_meshQuality} \&\ref{fig.16:wing_meshQuality2}. The orthogonality angle of iced and sinusoidal deformed mesh is compared to that of the mesh prior to deformation. The first location of interest is in the $x-z$ plane at $y = 0$ and is shown in Fig.~\ref{fig.15:wing_meshQuality}. This cut-plane provides a convenient view of the iced mesh. It shows there is a slight decrease in the mesh quality along the leading edge of the iced mesh. Despite this, the overall impact of the ice accretion on the mesh quality appears low. In this plane of view it is difficult to asses the impact of the sinusoidal deformation on the mesh quality. The second location of interest is in the $y-z$ plane at $x = c$ and is shown in Fig.~\ref{fig.16:wing_meshQuality2}. Likewise to the previous view, this exposes the marginal reduction in the mesh quality at the leading edge of the iced mesh. This cut-plane however provides a more favourable view of the sinusoidal deformation. Resultantly, it can be adjudged that the mesh quality of the sinusoidal deformation is also satisfactory. Overall it can be determined that both kind of deformation maintain a high mesh quality. 

\section{Conclusion}
\label{sec:conclusion}
Concurrent radial basis function mesh deformation numerical techniques are evaluated for their use within aircraft icing simulations. In this work state-of-the-art multi-level greedy surface point reduction and volume point reduction algorithms are used. To illustrate the efficiency and robustness of the approach during icing conditions a two-dimensional airfoil test case and a three-dimensional swept-wing test case are used. The data-reduction methods help to reduce the computational cost of the deformation process while maintaining high mesh quality. The method allows for the isolation of control points where the ice is present. The convergence history of the airfoil test case performs well with limited control points required at each level and subsequently ensures high performance. The convergence history of the much larger swept-wing test case also performs well and is able to achieve three levels of greedy surface point selection. The introduction of these techniques to aircraft icing simulations provides a significant improvement to the efficiency of standard radial basis function mesh deformation. Simultaneously, this work highlights the promise of radial basis function mesh deformation techniques and hopes to help provide a reasonable solution to the challenge that is accounting for the moving ice boundary. 

\section*{Acknowledgments}
The NITROS (Network for Innovative Training on ROtorcraft Safety) project has received funding from the European Union's Horizon 2020 research and innovation program under the Marie Sk\l{}odowska-Curie grant agreement No. 721920.

The ICE GENESIS project has received funding from the European Union's Horizon 2020 research and innovation program under grant agreement No. 824310.


\newpage
\section*{Figures}

\clearpage
\section*{Tables}

\begin{table}[htb]
\centering
\caption[Wendland radial basis functions with compact support]{Wendland radial basis functions with compact support where, $ \eta = \left ( \left \| \textbf{r} - \textbf{r}_{i} \right \| /R \right )$.}
\label{tab.1:rbf_functions}
\begin{tabular}{ll} \hline \hline
Name        & Basis Function \\ \hline
Wendland C0 & $ \phi(\eta) = \left ( 1 - \eta \right )^{2} $\\
Wendland C2 & $ \phi(\eta) = \left ( 1 - \eta \right )^{4} \left ( 4 \eta + 1 \right ) $\\
Wendland C4 & $ \phi(\eta) = \left ( 1 - \eta \right )^{6} ( \frac{35}{3} \eta^2 + 6 \eta + 1 ) $ \\
Wendland C6 & $ \phi(\eta) = \left ( 1 - \eta \right )^{8} ( 32 \eta^3 + 25 \eta^2 + 8 \eta + 1 ) $ \\ \hline \hline
\end{tabular}
\end{table}

\begin{landscape}

\begin{table}[h!]
\centering
\caption{NACA0012 icing conditions.}
\label{tab:2}
\begin{tabular}{lllllllll}
\hline \hline
\begin{tabular}[c]{@{}l@{}}Airfoil\\ {[}-{]}\end{tabular} & \begin{tabular}[c]{@{}l@{}}Accretion \\ Time\\ {[}s{]}\end{tabular} & \begin{tabular}[c]{@{}l@{}}Angle of \\ Attack\\ {[}$^\circ${]}\end{tabular} & \begin{tabular}[c]{@{}l@{}}Chord\\ {[}m{]}\end{tabular} & \begin{tabular}[c]{@{}l@{}} Pressure\\ {[}Pa{]}\end{tabular} & \begin{tabular}[c]{@{}l@{}}Airspeed\\ {[}m/s{]}\end{tabular} & \begin{tabular}[c]{@{}l@{}}Outside Air \\ Temperature\\ {[}K{]}\end{tabular} & \begin{tabular}[c]{@{}l@{}}Liquid Water \\ Content\\ {[}g/m$^{3}${]}\end{tabular} & \begin{tabular}[c]{@{}l@{}}Mean Volume \\ Diameter\\ {[}$\mu$m{]}\end{tabular} \\ \hline
NACA0012 & 120 & 0 & 0.3 & 90,700 & 129 & 260.55 & 0.5 & 20 \\ \hline \hline
\end{tabular}
\end{table}

\newpage

\begin{table}[htb!]
\centering
\caption{Swept wing icing conditions.}
\label{tab:3}
\begin{tabular}{lllllllll}
\hline \hline
\begin{tabular}[c]{@{}l@{}}Wing \\ Profile \\ {[}-{]}\end{tabular} & 
\begin{tabular}[c]{@{}l@{}}Sweep \\ Angle \\ {[}$^\circ${]}\end{tabular} &
\begin{tabular}[c]{@{}l@{}}Accretion \\ Time\\ {[}min{]}\end{tabular} & \begin{tabular}[c]{@{}l@{}}Angle of \\ Attack\\ {[}$^\circ${]}\end{tabular} & \begin{tabular}[c]{@{}l@{}}Chord\\ {[}m{]}\end{tabular} & 
\begin{tabular}[c]{@{}l@{}}Airspeed\\ {[}m/s{]}\end{tabular} & 
\begin{tabular}[c]{@{}l@{}}Total Air \\ Temperature\\ {[}K{]}\end{tabular} & \begin{tabular}[c]{@{}l@{}}Liquid Water \\ Content\\ {[}g/m$^{3}${]}\end{tabular} & \begin{tabular}[c]{@{}l@{}}Mean Volume \\ Diameter\\ {[}$\mu$m{]}\end{tabular} \\ \hline
NACA0012 & 45 & 19.8 & 0 & 0.914 & 51.44 & 260.55 & 0.57 & 44 \\ \hline \hline
\end{tabular}
\end{table}

\end{landscape}

\newpage 
\listoffigures

\newpage
\listoftables

\newpage
\section*{References}
\bibliography{mybibfile}

\end{document}